\documentclass[reqno]{amsart}

\usepackage{appendix}
\usepackage[french,english]{babel}
\usepackage[utf8]{inputenc}
\usepackage[T1]{fontenc}
\usepackage[hypertexnames=false]{hyperref}
\usepackage{cite}
\usepackage{comment}
\usepackage{palatino}
\usepackage[autostyle]{csquotes}
\usepackage{amsmath}
\usepackage{mathtools}
\usepackage{amsfonts,amssymb,amsmath,latexsym,mathrsfs,bbm,stmaryrd}
\usepackage[all]{xy}
\usepackage[usenames]{color}
\usepackage{epsfig}
\usepackage{graphicx}
\usepackage{subcaption}
\usepackage{float,enumerate}
\usepackage{amsthm}
\usepackage{thmtools}
\usepackage{thm-restate}
\usepackage{times}
\usepackage{tikz}
\usetikzlibrary{arrows.meta,decorations.markings}

\theoremstyle{plain}
\newtheorem{theorem}{Theorem}[section]

\newtheorem{proposition}[theorem]{Proposition}

\newtheorem{lemma}[theorem]{Lemma}
\newtheorem{corollary}[theorem]{Corollary}

\theoremstyle{definition}
\newtheorem{definition}[theorem]{Definition}
\newtheorem{example}[theorem]{Example}
\newtheorem{remark}[theorem]{Remark}

\numberwithin{equation}{section}

\def\E{{\mathbb E}}

\def\R{{\mathbb R}}

\def\Z{{\mathbb Z}}

\def\CA{{\mathcal A}}

\def\CC{{\mathcal C}}

\def\<{{\langle}}
\def\>{{\rangle}}

\DeclareMathOperator{\cov}{Cov}

\DeclareMathOperator{\tr}{tr}
\DeclareMathOperator{\Tr}{Tr}

\begin{document}

\title[Global Fluctuations of Gaussian Elliptic Matrices]{Global Fluctuations of Gaussian Elliptic Matrices}

\author {
Lingxuan Wu and Zhi Yin
}

\address{
\parbox{\linewidth}{Lingxuan Wu,
School of Mathematics and Statistics, 
Central South University,\\
Changsha, Hunan 410083, China.\\
\texttt{lingxuan.wu@csu.edu.cn}}
}

\address{
\parbox{\linewidth}{Zhi Yin,
School of Mathematics and Statistics, 
Central South University,\\
Changsha, Hunan 410083, China.\\
\texttt{hustyinzhi@163.com}}
}

\date{\today}
\maketitle

\begin{abstract}
We introduce a spoke-arc decomposition of non-crossing annular pair partitions $NC_2(p,q)$ that records spoke type and orientation, isolates spoke-level contributions, and factorizes the dependence on the ellipticity parameter $\gamma$ into a spoke factor and arc weights. This yields closed-form descriptions of the limiting covariance of Gaussian elliptic matrices. 
As a corollary, we show that the independent family of Gaussian elliptic random matrices are asymptotically second-order free.
\end{abstract}

\bigskip




\section{Introduction and main results}
\label{sec:intro}

\subsection{Elliptic random matrix}

The elliptic random matrix, initially introduced by Girko \cite{Girko1986,Girko1995,Girko19952,Girko1997,Girko2006,Girko20062} is a natural interpolation between the GUE random matrix and the Ginibre random matrix. Namely, we have the following definition.
\begin{definition}\label{def:elliptic-matrix}
A \emph{complex $N \times N$ Gaussian elliptic matrix} of parameter $\gamma \in [-1, 1]$ is a random matrix $X= (X_{ij})_{i,j=1}^N$ such that the following hold:
\begin{enumerate}[\rm (1)]
\item $\{X_{ii}: 1\leq i \leq N \} \cup \{(X_{ij}, X_{ji}): 1\leq i<j \leq N\}$ is a collection of independent random elements,
\item $\{(X_{ij}, X_{ji}): 1\leq i <j \leq N \}$ are i.i.d. Gaussian, centered, such that  
\begin{equation*}
\mathbb{E} X_{ij}^2 = \mathbb{E} X_{ji}^2 = \mathbb{E} X_{ij} \overline{X}_{ji} =0,\; \mathbb{E} \vert X_{ij} \vert^2 = \mathbb{E} \vert X_{ji} \vert^2 =1, \; \text{and} \; \mathbb{E} X_{ij} X_{ji} = \gamma;
\end{equation*}
\item $\{X_{ii}, 1\leq i \leq N \}$ are i.i.d. Gaussian, centered, such that $\mathbb{E} \vert X_{ii} \vert^2 =1$ and $\mathbb{E} X_{ii}^2 =\gamma$. 
\end{enumerate}
\end{definition}

We will also consider the real case, which is an interpolation between the GOE model and the real Ginibre model. 
\begin{definition}\label{def:real-elliptic-matrix}
A \emph{real Gaussian elliptic matrix} of parameter $\gamma\in[-1,1]$ is a random matrix
$X=(X_{ij})_{i,j=1}^N$ such that the following hold:
\begin{enumerate}[\rm (1)]
\item $\{X_{ii}: 1\le i\le N\}\ \cup\ \{(X_{ij},X_{ji}): 1\le i<j\le N\}$ is a collection of independent random elements;
\item $\{(X_{ij},X_{ji}): 1\le i<j\le N\}$ are i.i.d.\ centered real Gaussian vectors such that
\[
\mathbb{E}X_{ij}^2=\mathbb{E}X_{ji}^2=1 \; \text{and} \; 
\mathbb{E}X_{ij}X_{ji}=\gamma;
\]
\item $\{X_{ii}: 1\le i\le N\}$ are i.i.d.\ centered real Gaussian random variables such that
\[
\mathbb{E}X_{ii}^2=1+\gamma.
\]
\end{enumerate}
\end{definition}

At $\gamma=1$, the elliptic model reduces to the GUE or GOE model, at $\gamma=-1$ it becomes a skew-Hermitian or skew-symmetric Gaussian ensemble whose spectrum lies on the imaginary axis, while at $\gamma=0$ it coincides with Ginibre model, so the parameter $\gamma\in[-1,1]$ continuously deforms between normal and non-normal behavior. Beyond serving as a formal interpolation, elliptic ensembles also arise naturally in applications where correlations between real and imaginary parts destroy exact Hermitian symmetry.  Examples include the dynamics of random neural networks \cite{Crisanti1987,Sommers1988,Marti2018}, ecological stability \cite{Fyodorov2016}, complex networks \cite{Poley2024}, and non-Hermitian models in quantum physics \cite{Schomerus2017,Fyodorov1997}. In this sense, elliptic ensembles provide a minimal but realistic laboratory for non-normality, with the deformation parameter $\gamma$ measuring the degree of deviation from normality.

\subsection{Moments and the spectral of the elliptic model}

From the spectral perspective, the global distribution of the elliptic random matrix is by now fully understood: Girko's elliptic law, rigorously proved by Nguyen and O'Rourke \cite{NguyenORourke2015} (for the complex case) and by Naumov \cite{Naumov12} (for the real case), extending the circular law for Ginibre matrices \cite{TaoVu2010, Bordenave-Chafai-circular}. More precisely, Let $\lambda_1(X), \lambda_2(X), \ldots,\lambda_N(X)$ be the eigenvalues of $X$. The empirical spectral distribution (ESD) $\mu_{X}$ of $X$ is a random probability measure defined as
\begin{equation}\label{eqn:ESD-X}
	    \mu_{X}=\frac{1}{N}\sum_{k=1}^{N}\delta_{\lambda_k(X)}.
\end{equation}

\medskip

\noindent\emph{\bf Convention.}
Unless explicitly stated otherwise, we will work with the \emph{normalized} elliptic matrix and suppress the size index, writing simply
\begin{equation}\label{eq:normalized-elliptic}
X \equiv X_N \coloneqq \frac{1}{\sqrt{N}}(X_{ij})_{i,j=1}^N,
\end{equation}
where $(X_{ij})$ is a (complex or real) Gaussian elliptic matrix with parameter $\gamma\in[-1,1]$ given by Definition \ref{def:elliptic-matrix} or \ref{def:real-elliptic-matrix}.

\begin{theorem}[The elliptic law]
Let $X$ be the normalized Gaussian elliptic matrix given in \eqref{eq:normalized-elliptic}, then the limit (as $N \to \infty$) ESD of $X$ is the uniform measure supported in the ellipsoid defined as
 \begin{equation}\label{eqn:ellipsoid}
 \mathcal{E}_\gamma : = \left\{ (x, y) \in \mathbb{R}^2 : \frac{x^2}{(1+\gamma)^2} + \frac{y^2}{(1-\gamma)^2} \leq 1\right\}.
 \end{equation}
 \end{theorem}

On the other hand, the moments of the elliptic random matrix have a combinatoric description. For $n\in\mathbb{N}$, we denote $[n]:=\{1,\dots,n\}$. 
Denote $NC_2(n)$ by the set of all non-crossing partitions of $[n].$ Let $\tau:[n]\to\{1,*\}$  be a type word of length $n.$ We refer to Section \ref{sec:notation} for the definitions. Then we have the following theorem, which was proved by Adhikari and Bose in \cite{adhikari2019brown}.
\begin{theorem}\label{thm:moment}
Let $X$ be the normalized Gaussian elliptic matrix given in \eqref{eq:normalized-elliptic}. Given a type word $\tau$ of length $p,$ for any $\pi \in NC_2(2k),$ denote
$s(\pi;\tau)\coloneqq\#\{ \{x,y\} \in\pi:\ \tau(x)=\tau(y)\}.$
Then we have 
\begin{equation}
\lim_{N \to \infty} \mathbb{E} \tr \left(X^{\tau(1)} X^{\tau(2)} \cdots X^{\tau(p)} \right) = \begin{cases} \sum_{\pi \in NC_2(2k)} \gamma^{s(\pi; \tau)}, &  \text{if} \;  \; p = 2k,\\ 0, & \text{otherwise},\end{cases}
\end{equation}
where $\tr\coloneqq \Tr/N$ is the normalized trace. 
 \end{theorem}
An immediate result of the above theorem is that the independent family of elliptic random matrices are asymptotically free. Asymptotic freeness of random matrices was firstly discovered by Voiculescu \cite{Voiculescu1991, Voiculescu1998}. It provides a conceptual framework to understand the convergence of ESDs of large-dimensional random matrix ensembles. Hence, the above two results show that, at the first-order level, the picture for the elliptic ensemble is complete. 
 

\subsection{Main results}
In this paper, we introduce a combinatorial argument based on a decomposition of non-crossing annular pair partitions into spokes and arcs, which we call the spoke–arc decomposition (see Section~\ref{subsec:spoke-arc}). It yields closed-form covariance formulas for several canonical word families, and thus provides a quantitative description of fluctuations of elliptic ensembles across the entire ellipticity range. 

Now we are ready to present our main results, and the proofs are postponed to Section \ref{sec:proof}. In contrast to the first-order, there is a prominent difference between the complex and real cases for the results of the second-order. Let $\xi_1, \xi_2$ be two random variables, their covariance is given by 
\begin{equation*}
\kappa_2(\xi_1,\xi_2)=\cov (\xi_1,\xi_2)=\E(\xi_1 \xi_2)-\E(\xi_1)\E(\xi_2).
\end{equation*}

\begin{theorem}\label{thm:spoke-arc-closed}
Let $X$ be the normalized complex Gaussian elliptic matrix given in \eqref{eq:normalized-elliptic}.
For integers $p,q\ge1$, the limiting covariances of the three typical words are given as follows:
\begin{itemize}
\item Pure powers:
\begin{equation}\label{eq:FCpp}
\begin{split}
\lim_{N\to\infty}\kappa_2 \left(\Tr X^{p},\Tr X^{q} \right)
&=\sum_{a=1}^{\min\{p,q\}} a  \gamma^{\frac{p+q}{2}}\,\binom{p}{\tfrac{p-a}{2}} \binom{q}{\tfrac{q-a}{2}}.
  \end{split}
\end{equation}

\item Pure power vs.\ adjoint power:
\begin{equation}\label{eq:FCpa}
\begin{split}
\lim_{N\to\infty}\kappa_2\left(\Tr X^{p},\Tr (X^{\ast})^{q}\right)
&=\sum_{a=1}^{\min\{p,q\}} a \gamma^{\frac{p+q}{2}-a}\,
  \binom{p}{\tfrac{p-a}{2}}
  \binom{q}{\tfrac{p-a}{2}}.
  \end{split}
  \end{equation}
  
 \item Alternating:
 \begin{equation}\label{eq:FCall}
 \begin{split}
\lim_{N\to\infty}\kappa_2 &\left(\Tr(XX^{\ast})^{p},\Tr(XX^{\ast})^{q}\right)\\
& =\sum_{a=1}^{\min(p,q)}
a(\gamma^{2a}+1)
\binom{2p}{p-a}
\binom{2q}{q-a}.
   \end{split}
   \end{equation}
\end{itemize}
\end{theorem}

For $\gamma=1,$ $X$ reduces to the GUE model. It is known that (see \cite{MingoNica2004,Redelmeier2014,tutte1962census})
\begin{equation*}
\begin{split}
\lim_{N\to\infty}\kappa_2 \left(\Tr X^{p},\Tr X^{q} \right) &= \lim_{N\to\infty}\kappa_2\left(\Tr X^{p},\Tr (X^{\ast})^{q}\right)\\
& = |NC_2(p,q)| \\
& =\frac{1+(-1)^{p+q}}{2}\cdot \frac{2\,\lceil p/2\rceil\,\lceil q/2\rceil}{p+q}\,
\binom{p}{\lfloor p/2\rfloor}\binom{q}{\lfloor q/2\rfloor},
\end{split}
\end{equation*}
and 
\begin{equation*}
\lim_{N\to\infty}\kappa_2 \left(\Tr(XX^{\ast})^{p},\Tr(XX^{\ast})^{q}\right) = |NC_2(2p,2q)|, 
\end{equation*}
where $NC_2(p, q)$ is the set of all non-crossing annular pair partitions. Hence, by \eqref{eq:FCpp} we have the following 
interesting combinatoric identity. 
\begin{corollary}
Suppose that $p,q\ge 1$ have the same parity, we have 
\begin{equation}
|NC_2(p,q)|
=\sum_{a=1}^{\min\{p,q\}}\! a\,
\binom{p}{\tfrac{p-a}{2}} \binom{q}{\tfrac{q-a}{2}}.
\end{equation}
\end{corollary}

For $\gamma=0,$ now $X$ reduces to the Ginibre model. It is clear that 
\begin{equation*}
\lim_{N\to\infty}\kappa_2 \left(\Tr X^{p},\Tr X^{q} \right)  =0.
\end{equation*}
On the other hand, since only the term $a= (p+q)/2$ survives in the summation of \eqref{eq:FCpa}, we obtain 
\begin{equation*}
 \lim_{N\to\infty}\kappa_2\left(\Tr X^{p},\Tr (X^{\ast})^{q}\right) = \delta_{p, q}\cdot p. 
\end{equation*}
Combinatorially, this corresponds precisely to counting the non-crossing annular pair partitions in which the pairs are all mixed-type.
Also, for the alternating case, the limit of the covariance $\kappa_2 \left(\Tr(XX^{\ast})^{p}, \Tr(XX^{\ast})^{q}\right)$ reduces to counting such mixed-type partitions.  In particular, \cite{Redelmeier2014} explicitly describes the same combinatorial structures in the context of the real Ginibre ensemble, where the fluctuations are governed by non-crossing annular $\ast$-pairings---partitions that link untransposed and transposed terms.

\medskip

We now turn to the real case.
\begin{theorem}\label{thm:spoke-arc-closed-real}
Let $X$ be the normalized real Gaussian elliptic matrix given in \eqref{eq:normalized-elliptic}.
For integers $p,q\ge1$, the limiting covariances of the three typical words are given as follows:
\begin{itemize}
\item Pure powers and power vs.\ adjoint power:
\begin{equation}\label{eq:FCpp-real}
\begin{split}
\lim_{N\to\infty}\kappa_2\left(\Tr X^{p},\Tr X^{q}\right) &= \lim_{N\to\infty}\kappa_2\left(\Tr X^{p},\Tr (X^{\ast})^{q}\right)\\
=&\sum_{a=1}^{\min\{p,q\}} a\gamma^{\frac{p+q}{2}-a} (\gamma^{a}+1)
  \binom{p}{\tfrac{p-a}{2}} \binom{q}{\tfrac{q-a}{2}}.
  \end{split}
  \end{equation}
 
 \item Alternating:
 \begin{equation}\label{eq:FCall-real}
 \begin{split}
\lim_{N\to\infty}\kappa_2&\left(\Tr(XX^{\ast})^{p},\Tr(XX^{\ast})^{q}\right)\\
&=\sum_{a=1}^{\min\{p,q\}} 2a(\gamma^{2a}+1)
\binom{2p}{p-a}\binom{2q}{q-a}.
   \end{split}
\end{equation}
\end{itemize}
\end{theorem}

For $\gamma =0$ or $1,$ comparing with \eqref{eq:FCpp}-\eqref{eq:FCall}, it is clear that the limiting covariance of the real model is \emph{twice} that of the complex model. Except 
\begin{equation*}
\lim_{N\to\infty}\kappa_2\left(\Tr X^{p},\Tr X^{q}\right)=\delta_{p,q} \cdot p
\end{equation*}
for $\gamma=0,$ while the value is equal to $0$ for the complex case.

\medskip

In order to compute the above limiting covariance, we adopt a spoke-arc locality viewpoint: the dominant term of the covariance is determined by some admissible annular pair partitions.  
They are decomposed into local spoke/arc parts and the argument is carried out entirely at this local level. In principle, our method can be applied to the covariances of traces of more general words in $X$ and $X^*$. We refer to Section \ref{sec:spoke-arc} for more details about the spoke-arc decomposition. 

\medskip

As an application, our method can be used to study the asymptotic second-order freeness of Gaussian elliptic matrices. We refer to Section \ref{subsec:free} for the notion of (second-order) freeness. 
A key combinatorial precursor is the work of Mingo and Nica \cite{MingoNica2004}, which introduced and analyzed non-crossing annular permutations and partitions, laying the groundwork for second-order theory. Building on these ideas, a systematic theory of (complex) second-order freeness was developed in a trilogy of works \cite{MingoSpeicher2006, Mingo2007, Collins2007}, which established the framework (Part~I), extended it to unitary ensembles (Part~II), and further pushed the theory to higher orders via free cumulants (Part~III). Later, the notion of real second-order freeness was introduced by Redelmeier \cite{Redelmeier2014} and was subsequently developed in the context of Haar orthogonal matrices by Mingo and Popa \cite{MingoPopa2013}. 

It is worth emphasizing that at the second-order level the real and complex settings diverge markedly. In the complex case, only orientation-preserving non-crossing configurations contribute, and the definition of second-order freeness from \cite{MingoSpeicher2006} applies directly. In the real case, additional orientation-reversing contributions appear (already visible in GOE, real Wishart, and real Ginibre), which necessitate the modified notion of real second-order freeness \cite{Redelmeier2014}.

In this paper, we prove the following theorem. 

\begin{theorem}\label{thm:elliptic-r2}
For each size $N$, let $\{X^N_i\}_{i\in I}$ be a family of independent Gaussian elliptic matrices (with possibly different parameters). Then the family $\{X^N_i\}_{i\in I}$ is
second-order asymptotically free both in the complex and the real cases. 
\end{theorem}

The above theorem is a corollary of Mingo--\'Sniady--Speicher \cite[Corollary~3.16]{Mingo2007} (for the complex case) and Mingo--Popa \cite[Theorem~54]{MingoPopa2013} (for the real case).
For instance, in the real setting, the real Gaussian elliptic random matrix can be written as a linear combination of two independent random matrices: one is a real Ginibre Gaussian matrix and the another is a GOE matrix. Mingo--Popa showed that if two independent families each have a real second-order limit distribution and one family is invariant under orthogonal conjugation, then the two families are asymptotically real second-order free. Hence, it \emph{qualitatively} implies the real part of Theorem \ref{thm:elliptic-r2}.

However, our contribution here is \emph{quantitative}. We provide explicit formulas for the limiting covariances, and, combined with a detailed combinatorial analysis of annular non-crossing pair partitions, these formulas allow us to verify directly the asymptotic second-order freeness of independent Gaussian elliptic ensembles. In particular, in the real case, the contribution of the additional orientation-reversing non-crossing configurations can be clearly captured via our spoke-arc decomposition. 
We refer to Section \ref{subsec:second-free} for more details. 

\subsection{Organization of the Paper}
The remainder of the paper is organized as follows: 
Section~\ref{sec:notation} introduces necessary notions used throughout. 
Section~\ref{sec:spoke-arc} develops the spoke-arc decomposition and explains how it leads to factorization of $\gamma$-weights. 
Section~\ref{sec:proof} contains the proofs of the main theorems, combining combinatorial and analytic arguments.

\section{Preliminaries}
\label{sec:notation}

\subsection{Some combinatorics}
In this section, we will recall some necessary notions in combinatorics, which can be found in \cite{nica2006lectures, mingo2017free}. For $n\in\mathbb{N}$, we denote $[n]:=\{1,\dots,n\}$. 

\begin{definition}
A \emph{type word} of length $n$ is a map $\tau:[n]\to\{1,*\}$, where $\{1,*\}$ is the \emph{type alphabet}.
For $i\in[n]$, the symbol $\tau(i)\in\{1,*\}$ is called the $i$-th \emph{letter} (of type $1$ or $*$). 
The number of letters (the \emph{letter length}) is $|\tau|\coloneqq n$; by convention, the empty word has $|\tau|=0$.
Given such a $\tau$, we associate to it the matrix monomial
\[
X^{\tau(1)} X^{\tau(2)} \cdots X^{\tau(n)}.
\]
In this way, type words serve as a combinatorial encoding of monomials in $(X,X^\ast)$.

Moreover, if $\tau(i) =1$ or $\ast$ for any $i \in [n],$ then $\tau$ is called \emph{pure}. If $\tau(2j-1)=1$ and $\tau(2j)=\ast$ for $j=1, \ldots, n,$ then $\tau$ is called \emph{alternating}. 
\end{definition}

\begin{definition}
Let $\mathcal C$ be a finite set of colors. For each $c\in\mathcal{C}$ we fix a random matrix ensemble $X^{(c)}$. Distinct colors index distinct ensembles. A \emph{cluster} of color $c$ is a single-color matrix monomial
\[
Y = \prod_{i=1}^{n} \big(X^{(c)}\big)^{\tau(i)}, \quad \tau:[n]\to\{1,*\}.
\]
Its letter length is $|Y|:=|\tau|=n$ (by convention $|I|=0$).
\end{definition}

\begin{remark}
On annular partition diagrams, a cluster $Y$ of color $c$ corresponds to a letter interval $I$ labeled by the
alphabet $\{X^{(c)},(X^{(c)})^{*}\}$ according to $\tau$; we write $c=c(I)=c(Y)$ and $\tau(I)=\tau$.
\end{remark}

\begin{definition}
A \emph{partition} of the finite set $[n]$ is a collection 
\[
\pi=\{V_1,V_2,\dots,V_r\}
\]
of nonempty disjoint subsets $V_i \subseteq [n]$ such that $V_1 \cup V_2 \cup \cdots \cup V_r=[n]$.  
The subsets $V_i$ are called the \emph{blocks} of $\pi$.  
We denote by $\mathcal P(n)$ the set of all partitions of $[n]$. A partition is called a \emph{pair partition} if every block contains exactly two elements. We denote the set of pair partitions on $[n]$ by $\mathcal P_2(n)$. If $n$ is odd, then $\mathcal P_2(n)=\varnothing$.
\end{definition}

We will denote the group permutations on $[n]$ by $S_n$. We always use the cycle notation for group permutations. Unless otherwise specified, we identify a partition $\pi \in \mathcal P(n)$ with a permutation $\sigma \in S_n$ by writing each block as a cycle, up to cyclic rotation within each cycle. Whenever a construction depends on the cyclic order, we will state the choice explicitly. We note that a pair partition in $\mathcal P_2(n)$ uniquely defines a permutation in $S_n$. We denote the number of cycles of a permutation $\pi$ by $\#(\pi)$.

\begin{definition}\label{def:non-crossing}
A partition $\pi \in \mathcal P(n)$ is called \emph{non-crossing} if it contains no crossing: there do not exist indices
\[
1\le a<b<c<d\le n
\]
and two distinct blocks $B,C\in\pi$ such that $a,c\in B$ and $b,d\in C$.
Equivalently, place the points $1,\dots,n$ in clockwise order on a circle and, for each block of $\pi$, connect its elements by chords inside the circle.
Then $\pi$ is non-crossing if and only if these chords can be drawn so that no two chords intersect.
We denote by $NC(n)$ the set of all non-crossing partitions of $[n]$. The set of \emph{non-crossing pair partitions} is denoted by $NC_2(n)\subseteq\mathcal P_2(n)$.

\end{definition}

It is well-known that $|NC(n)|=|NC_2(2n)|=C_n$, which is the $n$-th Catalan number with generating function
$\mathcal C(z)=\sum_{n\ge 0} C_n z^n=\frac{1-\sqrt{1-4z}}{2z}$.

\begin{figure}[H]
\centering
\begin{tikzpicture}[line cap=round,line join=round,thick,scale=.8]
\tikzset{
  directed/.style = {postaction={decorate},
    decoration={markings, mark=at position 0.55 with {\arrow{Stealth[length=2mm]}}}},
  crosscycle/.style = {line width=0.8pt, directed},
}

\def\R{2.4}   
\def\r{1}     

\draw (0,0) circle (\R);
\draw (0,0) circle (\r);

\node at (30:\R+0.30)  {$10$};
\node at (90:\R+0.30)  {$5$};
\node at (150:\R+0.30) {$6$};
\node at (210:\R+0.30) {$7$};
\node at (-30:\R+0.30) {$9$};
\node at (-90:\R+0.30) {$8$};

\node at (90:\r-0.30)  {$1$};
\node at (0:\r-0.30)   {$2$};
\node at (-90:\r-0.30) {$3$};
\node at (180:\r-0.30) {$4$};

\coordinate (10) at (30:\R);
\coordinate (5)  at (90:\R);
\coordinate (6)  at (150:\R);
\coordinate (7)  at (210:\R);
\coordinate (8)  at (-90:\R);
\coordinate (9)  at (-30:\R);
\coordinate (1) at (90:\r);
\coordinate (2) at (0:\r);
\coordinate (3) at (-90:\r);
\coordinate (4) at (180:\r);

\draw[crosscycle] (1) to[out=30,in=-60,looseness=1] (5);
\draw[crosscycle] (5) to[out=-30,in=180,looseness=1] (10); 
\draw[crosscycle] (10) to[out=180,in=30,looseness=1] (2); 
\draw[crosscycle] (2) to[out=60,in=30,looseness=1.4] (1);

\draw[crosscycle] (4) to[out=-110,in=200,looseness=1.2] (3);
\draw[crosscycle] (3) to[out=-150,in=-120,looseness=1.8] (4);

\draw[crosscycle] (6) to[out=-90,in=210,looseness=1.4] (9);
\draw[crosscycle] (9) to[out=220,in=20,looseness=1.0] (8);
\draw[crosscycle] (8) to[out=140,in=-20,looseness=1.0] (7);
\draw[crosscycle] (7) to[out=100,in=-100,looseness=1.0] (6);

\end{tikzpicture}
\caption{Example of a non-crossing annular permutation highlighting cycle orientation. The instance shown is $\pi=(1\,5\,10\,2)(3\,4)(6\,9\,8\,7)\in S_{NC}(4,6)$.}
\label{fig:annular-example-general}
\end{figure}
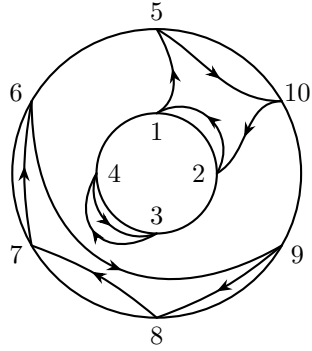

\begin{definition}[Non-crossing annular permutation \cite{MingoNica2004,mingo2017free}]\label{def:ncannular}
Fix $p,q \in \mathbb{N}$. Place the points $1,\dots,p$ on an inner circle clockwise and the points $p+1,\dots,p+q$ on a concentric outer circle counterclockwise of a annulus. For a permutation $\pi\in S_{p+q}$, we draw the cycles of $\pi$ between the circles by connecting consecutive elements of each cycle with chords. 
We say $\pi$ is a \emph{non-crossing annular permutation} if
\begin{itemize}
	\item the cycles do not cross,
	\item at least one cycle connects the two circles, and
	\item each cycle encloses a region between the circles homeomorphic to the disc with boundary oriented clockwise.
\end{itemize}
The set of all such permutations is denoted $S_{NC}(p,q)$. The subset consisting of non-crossing pair partitions is denoted by $NC_2(p,q)$. Pairs joining one outer point to one inner point are called \emph{spokes}. Pairs joining the points on the same circle are called \emph{arc pairs}. Throughout, we assume $p$ and $q$ have the same parity. See Fig \ref{fig:annular-example-general} for an example.
\end{definition}

\begin{remark}\label{rmk:eqdef-ncannular}
	 In \cite[Section~6]{MingoNica2004} it was shown that under the connectedness condition (i.e., at least one cycle connects the two circles), this definition was equivalent to the algebraic condition $\#(\pi)+\#(\pi^{-1}\rho)=p+q$, where $\rho=(1\cdots p)(p+1\cdots p+q)$.
\end{remark}


\subsection{Non-commutative probability spaces and freeness}\label{subsec:free}
In order to describe the limiting behavior of random matrices in the large-dimensional regime, we recall the basic notions of free probability, and we refer to \cite{MingoSpeicher2006,mingo2017free,nica2006lectures,Redelmeier2014} for more details.
\begin{definition}\label{def:ncps}
A \emph{non-commutative probability space} is a pair $(\mathcal A,\varphi)$ where
$\mathcal A$ is a unital algebra and $\varphi:\mathcal A\to\mathbb C$ is a unital linear functional satisfying
  $\varphi(1)=1$.
Elements of $\mathcal A$ are called (non-commutative) random variables. If, in addition, $\varphi$ is \emph{tracial} ($\varphi(ab)=\varphi(ba)$ for all $a,b\in\mathcal A$),
we say $(\mathcal A,\varphi)$ is \emph{tracial}.

A \emph{second-order non-commutative probability space} $(\CA, \varphi_1,\varphi_2)$ is a non-commutative probability space $(\CA,\varphi_1)$ equipped with a bilinear functional $\varphi_2$ that is tracial in each variable, symmetric in its two arguments and such that $\varphi_2(1,a)=\varphi_2(a, 1)=0$ for all $a\in \CA$.
\end{definition}

\begin{remark}
	We say a random variable $a$ is \emph{centered} if $\varphi_1(a)=0$. If $X,Y$ are $N\times N$ random matrices, we let $\varphi_1(X)=\E\tr(X)$ and $\varphi_2(X,Y)=\cov(\Tr(X),\Tr(Y))$, where $\tr\coloneqq \Tr/N$ is the normalized trace. 
\end{remark}

\begin{definition}
Let $\mathcal C$ be a finite color set, and for each $c\in\mathcal C,$ let $\mathcal A^{(c)}$ denote the subalgebra indexed by the color $c$.
Given a finite sequence of elements $a_1,\dots,a_n$ with $a_r\in\mathcal A^{(c(r))}$ for some color map $c:[n]\to\mathcal C$, we say that $(a_1,\dots,a_n)$ is \emph{alternating} if
$
c(r)\neq c(r+1)$ for all $r=1,\dots,n-1$. The sequence $(a_1,\dots,a_n)$ is called \emph{cyclically alternating} if it is alternating and in addition
$c(n)\neq c(1)$. 
\end{definition}

The notion of freeness plays the same role in free probability as independence does in classical probability.  It provides a usable ``calculus'' for combining non-commutative variables: under freeness, the distribution of $a+b$ (resp. $ab$) is determined by the individual laws of $a$ and $b$ via the free additive (resp. multiplicative) operation---mirroring how classical convolution controls sums of independent variables. 
 \begin{definition}\label{def:free}
Let $(\mathcal A,\varphi)$ be a non-commutative probability space and let
$\{\mathcal A^{(c)}\}_{c\in\mathcal C}$ be the subalgebras of $\mathcal A$ indexed by a finite color set $\mathcal C$.
We say that the family $\{\mathcal A^{(c)}\}_{c\in\mathcal C}$ is \emph{free} if for every $n\ge1$ and every
\emph{alternating} sequence $(a_1,\dots,a_n)$ with $a_r\in\mathcal A^{(c(r))}$ and $\varphi(a_r)=0$ for all $r$, one has
\[
\varphi(a_1, \ldots,a_n)=0.
\]
Families of random variables $\{a_i\}_{i}$ in $\mathcal A$ are said to be \emph{free} if the subalgebras they generate form a free family.
\end{definition}

Freeness explains first-order relations among moments. To analyze correlations beyond expectation, we pass to a second functional $\varphi_2$ and the corresponding notion of second-order freeness, which elevates the independence-like calculus from moments to covariances.

\begin{definition}[Second-order freeness]\label{def:second-order-freeness}
Let $(\mathcal A,\varphi_1,\varphi_2)$ be a second-order non-commutative probability space, and let $\{\mathcal A^{(c)}\}_{c\in\mathcal C}$ be subalgebras of $\mathcal A$ indexed by a finite color set $\mathcal C$.
We say that $\{\mathcal A^{(c)}\}_{c\in\mathcal C}$ is \emph{second-order free} if:
\begin{itemize}
  \item the family is free with respect to $\varphi_1$.
  \item for any $p,q\ge1$ and any \emph{centered, cyclically alternating} sequences
  $(a_1,\dots,a_p)$ and $(b_1,\dots,b_q)$ with $a_i\in\mathcal A^{(c(i))}$,
  $b_j\in\mathcal A^{(d(j))}$, and $\varphi_1(a_i)=\varphi_1(b_j)=0$ for all $i,j$, one has
  \[
    \varphi_2(a_1\cdots a_p, b_1\cdots b_q)
    =
    \delta_{p,q}\sum_{k=1}^{p}\;\prod_{i=1}^{p}\varphi_1\big(a_i b_{\,k-i}\big),
  \]
  where the indices are read modulo $p$.
\end{itemize}
\end{definition}

\begin{definition}[Real second-order freeness]\label{def:real-second-order-freeness}
Let $(\mathcal A,\varphi_1,\varphi_2)$ be a second-order non-commutative probability space equipped with an involution $a \mapsto a^t$ reversing the order of multiplication, and let $\{\mathcal A^{(c)}\}_{c\in\mathcal C}$ be subalgebras of $\mathcal A$ indexed by a finite color set $\mathcal C$.
We say that $\{\mathcal A^{(c)}\}_{c\in\mathcal C}$ is \emph{real second-order free} if:
\begin{itemize}
  \item the family is free with respect to $\varphi_1$.
  \item for any $p,q\ge1$ and any \emph{centered, cyclically alternating} sequences
  $(a_1,\dots,a_p)$ and $(b_1,\dots,b_q)$ with $a_i\in\mathcal A^{(c(i))}$,
  $b_j\in\mathcal A^{(d(j))}$, and $\varphi_1(a_i)=\varphi_1(b_j)=0$ for all $i,j$, one has
  \[
    \varphi_2(a_1\cdots a_p,\; b_1\cdots b_q)
    =
    \delta_{p,q}\sum_{k=1}^{p}\Big(\prod_{i=1}^{p}\varphi_1\big(a_i b_{\,k-i}\big)+\prod_{i=1}^{p}\varphi_1\big(a_i b^t_{\,k+i}\big)\Big),
  \]
  where the indices are read modulo $p$.
\end{itemize}
\end{definition}

Let $\xi_1,\dots,\xi_n$ be scalar random variables. The \emph{$r$-th classical cumulant}
$\kappa_r$ is the $r$-linear functional determined by the moment-cumulant formula
\[
\E\Big[\prod_{j=1}^n \xi_j\Big]
=\sum_{\pi\in \mathcal P(n)}\ \prod_{B\in\pi}\ \kappa_{|B|}\big((\xi_j)_{j\in B}\big),
\]
where $|B|$ is the number of elements in the block $B$.
Equivalently, if 
\(K(t_1,\dots,t_n):=\log \E\big[\exp\big(\sum_{j=1}^{n} t_j\,\xi_j\big)\big]\)
is the cumulant generating function, then
\[
\kappa_r(\xi_{i_1},\dots,\xi_{i_r})
=\left.\frac{\partial^r}{\partial t_{i_1}\cdots\partial t_{i_r}}\,K(t_1,\dots,t_n)\right|_{t=0}.
\]
In particular, $\kappa_1(\xi)=\E(\xi)$ and $\kappa_2(\xi_1,\xi_2)=\cov (\xi_1,\xi_2)=\E(\xi_1 \xi_2)-\E(\xi_1)\E(\xi_2)$.

Although formulated abstractly, these notions are realized by high-dimensional random matrices, where $\varphi_1$ and $\varphi_2$ arise from normalized traces and trace covariances. In the large-$N$ limit, their convergence leads to first- and second-order limit distributions. We now make this precise and introduce (real) asymptotic second-order freeness.
\begin{definition}[Second-order limit distribution]\label{def:so-limit}
Let $\{X_i^{(N)}\}_{i\in I}$ be a family of random matrices. Denote by $\kappa_r$ the classical $r$-th cumulant. We say that $\{X_i^{(N)}\}_{i\in I}$ has a \emph{second-order limit distribution} if there exists a second-order non-commutative probability space $(\CA,\varphi_1,\varphi_2)$  with random variables $\{a_i\}_{i\in I}$ equipped with an involution $a\mapsto a^{*}$ on $\mathcal A$ such that
\begin{itemize}
  \item for any non-commutative $\ast$-polynomial $P\in \mathbb C\langle x_i, x^\ast_i\mid i\in I\rangle$,
  \[
  \lim_{N\to\infty}\E\big[\tr\big(P(X^{(N)})\big)\big] = \varphi_1(P(a)).
  \]
  \item for any non-commutative $\ast$-polynomials $P,Q \in \mathbb C\langle x_i, x^\ast_i \mid i\in I\rangle$,
  \[
  \lim_{N\to\infty}\kappa_2\Big(\Tr\big(P(X^{(N)})\big),\Tr\big(Q(X^{(N)})\big)\Big) 
  = \varphi_2(P(a),Q(a)).
  \]
  \item for every $r\ge 3$ and all $P_1,\ldots,P_r \in \mathbb C\langle x_i, x^\ast_i\mid i\in I\rangle$,  
  \[
  \lim_{N\to\infty}\kappa_r\Big(\Tr\big(P_1(X^{(N)})\big),\dots,\Tr\big(P_r(X^{(N)})\big)\Big)=0.
  \]
 \end{itemize}
 Here $P(X^{(N)})$ (resp.\ $P(a)$) denotes evaluation of the non-commutative $\ast$-polynomial $P$ under the substitution $x_i\mapsto X_i^{(N)}$, $x_i^{*}\mapsto X_i^{(N)*}$ (resp.\ $x_i\mapsto a_i$, $x_i^{*}\mapsto a_i^{\ast}$).
 \end{definition}

\begin{definition}[Asymptotic second-order freeness]\label{def:SO-asymp-free}
For each $N\in\mathbb N$, let $\{X^{(c)}_N\}_{c\in\CC}$ be a colored family of random matrices. We say that $\{X^{(c)}_N\}_{c\in\CC}$ are \emph{asymptotically second-order free} if they have a second-order limit distribution and the algebras generated by the corresponding limiting random variables are second-order free; equivalently, for any $p,q\geq 2$ and any centered, cyclically alternating sequences $(A_1,\dots,A_p)$ and $(B_1,\dots,B_q)$ with
\[
A_i\in \operatorname{alg}\langle X^{c(i)}_N,X^{c(i)\ast}_N\rangle,\quad
B_j\in \operatorname{alg}\langle X^{d(j)}_N,X^{d(j)\ast}_N\rangle,
\]
and $\mathbb E\tr(A_i)=\mathbb E\tr(B_j)=0$, one has
\[
\lim_{N\to\infty}\cov\Big(\Tr(A_1\cdots A_p),\Tr(B_1\cdots B_q)\Big)
=
\delta_{p,q}\sum_{k=1}^{p} \prod_{i=1}^{p} \lim_{N\to\infty}\E\tr\big(A_i B_{k-i}\big).
\]
\end{definition}

\begin{definition}[Asymptotic real second-order freeness]\label{def:SO-asymp-real-free}
For each $N\in\mathbb N$, let $\{X^{(c)}_N\}_{c\in\CC}$ be a colored family of random matrices. We say that $\{X^{(c)}_N\}_{c\in\CC}$ are \emph{asymptotically real second-order free} if they have a second-order limit distribution and the algebras generated by the corresponding limiting random variables are real second-order free; equivalently, for any $p,q\geq 2$ and any centered, cyclically alternating sequences $(A_1,\dots,A_p)$ and $(B_1,\dots,B_q)$ with
\[
A_i\in \operatorname{alg}\langle X^{c(i)}_N,X^{c(i)\ast}_N\rangle,\qquad
B_j\in \operatorname{alg}\langle X^{d(j)}_N,X^{d(j)\ast}_N\rangle,
\]
and $\mathbb E\tr(A_i)=\mathbb E\tr(B_j)=0$, one has
\begin{align}\label{eq:second-order-free-defin}
	&\lim_{N\to\infty}\cov\Big(\Tr(A_1\cdots A_p),\Tr(B_1\cdots B_q)\Big)\\ \notag
&\quad=
\delta_{p,q}\sum_{k=1}^{p} \Big( \prod_{i=1}^{p} \lim_{N\to\infty}\E\tr\big(A_i B_{k-i}\big)+\prod_{i=1}^{p} \lim_{N\to\infty}\E\tr\big(A_i B^T_{k+i}\big)\Big).
\end{align}
\end{definition}

\section{Spoke-arc decomposition and limiting covariances}\label{sec:spoke-arc}

\subsection{Spoke-arc decomposition}\label{subsec:spoke-arc}

We use the notations in \cite{mingo2017free}. For $p, q \in \mathbb{N},$ place the points $1,\dots,p$ on an inner circle clockwise and the points $p+1,\dots,p+q$ on a concentric outer circle counterclockwise of a annulus. 
Indices of points on the two circles are taken modulo $p$ and $q$, respectively. Let $\rho=(1, 2, \ldots, p)\,(p+1, p+2, \ldots,p+q)\in S_{p+q}$ be the trace permutation along these orientation.

First we isolate the basic picture behind non-crossing annular pair partitions on two circles. Any $\pi\in NC_2(p,q)$ can be read as a family of spokes (pairs joining the two circles) together with arc pairs (pairs that stay on one circle). Removing the spoke endpoints cuts each circle into consecutive open arcs; on every such arc the remaining points are matched by a non-crossing pair partition that is independent of the other arcs. The data that matter are:
(i) the number $a$ of spokes and the cyclic lists of their endpoints on the two circles, and (ii) the lengths of the induced arcs, i.e., the number of points between consecutive spoke endpoints. Up to a simultaneous cyclic relabeling of the spokes, these data encode $\pi$ uniquely. This ``spoke-arc'' decomposition is purely combinatorial and will be the mechanism behind the factorization of weights in the sequel.

\begin{definition}[Spoke]\label{def:spoke}
Let $\pi\in NC_2(p,q)$.  
A \emph{spoke} of $\pi$ is a pair $\{u,v\}\in\pi$ with $u\in [p]$ on the inner circle and $v\in p+[q]$ on the outer circle.  
If $\pi$ has $a$ spokes, we list their inner endpoints and outer endpoints on two circles along the clockwise order:
\[
U=(u_1,\dots,u_a)\subset [p],\qquad V=(v_1,\dots,v_a)\subset p+[q].
\]
The spokes connect $u_r$ with $v_r$ for $r=1,\dots,a$, with indices taken modulo $a$.  
\end{definition}

\begin{definition}[Arc]\label{def:arc}
Suppose that $\pi\in NC_2(p,q)$ have $a$ spokes with ordered endpoint lists 
$U=(u_1,\dots,u_a)\subset[p]$ and $V=(v_1,\dots,v_a)\subset p+[q]$.
Removing all spoke endpoints cuts each circle into $a$ open arcs.
For each $r=1,\dots,a$, let $I_r\subset[p]$ be the set of inner \emph{arc points} strictly between $u_r$ and $u_{r+1}$ on $[p]$, and let 
$O_r\subset p+[q]$ be the set of outer arc points strictly between $v_r$ and $v_{r+1}$ on $p+[q]$.
Write $\iota_r\coloneqq|I_r|$ and $o_r\coloneqq|O_r|$, so that
\[
\iota_1+\cdots+\iota_a=q-a,\quad o_1+\cdots+o_a=p-a,\quad 
\iota_r,o_r\in 2\mathbb N_0.
\]
An \emph{arc pair partition} on the $r$-th inner (resp.\ outer) arc means a non-crossing pair partition 
\[
\pi_r^{\boldsymbol{\iota}}\in NC_2(I_r)\quad(\text{resp.\ }\pi_r^{\mathbf o}\in NC_2(O_r)).
\]
\end{definition}

Thus $\pi$ is encoded by the data
\[
(a,\mathbf o,\boldsymbol\iota,U,V;\ \{\pi^{\mathbf o}_r\}_{r=1}^a,\ \{\pi^{\boldsymbol{\iota}}_r\}_{r=1}^a).
\]
The above tuple of data is called a \emph{labeled spoke-arc configuration}. It is also worth emphasizing that if we perform a simultaneous cyclic relabeling of spoke endpoints by $s\in\{0,\dots,a-1\}$, i.e.,
\[
(U,V)\ \longmapsto\ \big((u_{1+s},\dots,u_{a+s}),\ (v_{1+s},\dots,v_{a+s})\big),
\]
the resulting data describe the same pair partition $\pi$. With the above definitions in hand, any $\pi\in NC_2(p,q)$ admits the following decomposition.

\begin{definition}[Spoke-arc decomposition]\label{def:spoke-arc-decomposition}
The set $NC_2(p,q)$ of non-crossing annular pair partitions is parameterized by the bijection
\begin{align}
	\notag &NC_2(p,q)\cong\ \\
	&
\bigsqcup_{1\le a\le \min\{p,q\}}
\bigsqcup_{\substack{o_1+\cdots+o_a=p-a\\ \iota_1+\cdots+\iota_a=q-a}}
\Big(\mathcal U_{a,\mathbf o}\times\mathcal V_{a,\boldsymbol\iota}\times
\prod_{r=1}^a NC_2(o_r)\times\prod_{s=1}^a NC_2(\iota_s)\Big)\Big/\mathbb Z_a,
\end{align}
where $\mathcal U_{a,\mathbf o}$ (resp.\ $\mathcal V_{a,\boldsymbol\iota}$) is the set of ordered endpoint lists $U$ on $[p]$ (resp.\ $V$ on $p+[q]$) realizing the outer (resp.\ inner) arc lengths $\mathbf o$ (resp.\ $\boldsymbol \iota$).
We call this bijection the \emph{spoke-arc decomposition} of non-crossing annular pair partitions.
\end{definition}

\begin{remark}
Note that a simultaneous cyclic relabelling of both lists by $s\in\mathbb{Z}_a$ produces the same non-crossing annular pair partition $\pi\in NC_2(p,q)$. This is the reason why the bijection is $a$-to-$1$ before quotienting by $\mathbb Z_a$.
\end{remark}


\begin{figure}
\centering
\begin{tikzpicture}[line cap=round,line join=round,thick,scale=.8]

\def\R{2.4}   
\def\r{1}  

\draw (0,0) circle (\R);
\draw (0,0) circle (\r);
\node at (30:\R+0.30)  {$10$};
\node at (90:\R+0.30)  {$5$};
\node at (150:\R+0.30) {$6$};
\node at (210:\R+0.30) {$7$};
\node at (-30:\R+0.30) {$9$};
\node at (-90:\R+0.30) {$8$};
\node at (90:\r-0.30) {$1$};
\node at (0:\r-0.30) {$2$};
\node at (-90:\r-0.30) {$3$};
\node at (180:\r-0.30) {$4$};

\coordinate (10)  at (30:\R);
\coordinate (5)  at (90:\R);
\coordinate (6)  at (150:\R);
\coordinate (7)  at (210:\R);
\coordinate (8) at (-90:\R);    
\coordinate (9) at (-30:\R);   

\coordinate (1) at (90:\r);    
\coordinate (2) at (0:\r);
\coordinate (3) at (-90:\r);
\coordinate (4) at (180:\r);

\draw (1) to[out=30,in=210,looseness=1.2] (5);
\draw (3) to[out=-120,in=210,looseness=1.5] (4);
\draw (2) to[out=0,in=180,looseness=1.1] (10);
\draw (7) to[out=0,in=90,looseness=0.4] (8);
\draw (6) to[out=-90,in=-160,looseness=1.4] (9);
\end{tikzpicture}
\caption{An non-crossing annular pair partition $\pi\in NC_2(4,6)$ illustrating the spoke-arc decomposition.}
  \label{fig:annular-example}
\end{figure}
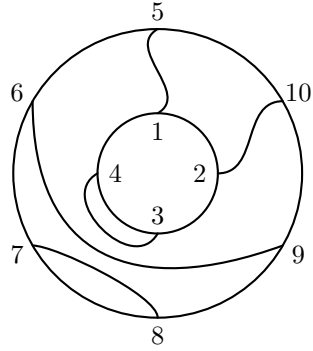


\begin{example}
Consider the partition
\[
\pi=\{\{1,5\},\{2,10\},\{3,4\},\{6,9\},\{7,8\}\}\in NC_2(4,6),
\]
depicted in Figure~\ref{fig:annular-example}. 

\noindent\emph{Spokes.} 
The cross-circle pairs are $\{1,5\}$ and $\{2,10\}$. Thus $a=2$, with ordered spoke endpoint lists
\[
U=(1,2)\subset [4],\quad V=(5,10)\subset 4+[6].
\]

\noindent\emph{Arcs.}
Removing these spoke endpoints cuts each circle into two arcs:
\begin{itemize}
\item On the inner circle $[4]$, the arc between $u_1=1$ and $u_2=2$ contains no points ($\iota_1=0$), while the arc between $u_2=2$ and $u_1=1$ contains $\{3,4\}$ ($\iota_2=2$).
\item On the outer circle $4+[6]$, the arc between $v_1=5$ and $v_2=10$ contains $\{6,7,8,9\}$ ($o_1=4$), while the arc between $v_2=10$ and $v_1=5$ is empty ($o_2=0$).
\end{itemize}

\noindent\emph{Arc pair partitions.}
The inner arc $\{3,4\}$ is matched by the pair $\{3,4\}$, giving $\pi^\iota_2=\{\{3,4\}\}$ and $\pi^\iota_1=\varnothing$.
On the outer arc $\{6,7,8,9\}$, the pairs are $\{6,9\}$ and $\{7,8\}$, giving $\pi^o_1=\{\{6,9\},\{7,8\}\}$ and $\pi^o_2=\varnothing$.

\noindent\emph{Cyclic relabeling.}
Since $a=2$, a cyclic shift by $1$ changes the endpoint lists to
\[
(U,V)=((1,2),(5,10))\ \longmapsto\ ((2,1),(10,5)),
\]
which swaps the roles of the two spokes. The resulting labeled configuration encodes the same partition $\pi$, showing explicitly how the $\mathbb Z_2$-action identifies different labelings with the same underlying diagram in Figure~\ref{fig:annular-example}.
\end{example}

\subsection{Semi-closed limiting covariance}\label{Semi-closed-lim-cov}
Recall
$\rho=(1,2, \ldots, p)\,(p+1, p+2, \ldots, p+q)\in S_{p+q}$.
Fix $p,q\ge1$ and a type word $\tau: [p+q]  \to \{1,\ast\}$.
Let $X$ be the normalized (complex or real) Gaussian elliptic matrix given in \eqref{eq:normalized-elliptic}. For any $t\in[p+q]$, we set
\begin{equation}\label{eq:W}
W_1\coloneqq \Tr\left(\prod_{t=1}^{p} X^{\tau(t)}\right),\quad
W_2\coloneqq \Tr\left(\prod_{t=p+1}^{p+q} X^{\tau(t)}\right).
\end{equation}
Then
\begin{equation}
W_1W_2=\sum_{i_1,\dots,i_{p+q}=1}^N\ \prod_{t=1}^{p+q}\ \big(X^{\tau(t)}\big)_{\,i_t\,i_{\rho(t)}}.
\end{equation}

\medskip
\subsubsection{Complex case}
\begin{proposition}\label{prop:semi-closed-cov}
Let $X$ be the normalized complex Gaussian elliptic matrix given in \eqref{eq:normalized-elliptic}, and let $W_1, W_2$ be given in \eqref{eq:W}. Then
\begin{equation}\label{eq:lim-cov-nc2-brief}
\lim_{N\to\infty}\kappa_2(W_1,W_2)=\sum_{\pi\in NC_2(p,q)} \gamma^{\,s(\pi;\tau)},
\end{equation}
where $s(\pi;\tau):=\#\{ \{x,y\} \in\pi:\ \tau(x)=\tau(y)\}$ counts the same-type pairs of $\pi$.
\end{proposition}

\begin{remark}
Formula \eqref{eq:lim-cov-nc2-brief} is \emph{semi-closed} in the sense that it gives an exact finite sum over the combinatorial objects $NC_2(p,q)$ with explicit monomial weights $\gamma^{\,s(\pi;\tau)}$, but it does not evaluate the partition count or the arc contributions further. In particular, the local arc contributions are not yet compressed. By contrast, the \emph{closed} forms in Theorem~\ref{thm:spoke-arc-closed} reorganize the sum by a single integer parameter and evaluate the arc part in closed form, producing a short sum of length $\min\{p,q\}$. Hence the dependence on $(p,q,\gamma)$ becomes explicit and numerically efficient.
\end{remark}

The proof of Proposition~\ref{prop:semi-closed-cov} is standard (see \cite{MingoSpeicher2006, Mingo2007}); we write down the proof for the completeness. The main technical tool is the following Wick's formula. A proof of Wick's formula can be found in \cite[p.~164]{Lando2004-dc}.

\begin{lemma}[Wick's formula]\label{lem:wick}
Let $(\xi_i)_{i\in I}$ be a centered jointly complex Gaussian family. 
For any $n\ge1$, indices $i_1,\dots,i_m\in I$, and a type word $\tau:[n]\to\{1,*\},$ we have
\[
\E\left[\prod_{u=1}^n \xi_{i_u}^{\,\tau(u)}\right]
=\begin{cases}
0,& n\ \text{odd},\\[2mm]
\displaystyle \sum_{\pi\in\mathcal P_2(n)} \prod_{\{x,y\}\in\pi}\ \E\left[\,\xi_{i_x}^{\,\tau(x)}\,\xi_{i_y}^{\,\tau(y)}\right],& n\ \text{even}.
\end{cases}
\]
\end{lemma}

We introduce the weights
\begin{equation}
w(1,1)=w(*,*)=\gamma,\quad w(1,*)=w(*,1)=1.
\end{equation}
Then for any $\alpha,\beta\in\{1,*\}$, by Definition \ref{def:elliptic-matrix}, the entrywise covariance reads
\begin{equation}\label{eq:elliptic-entry-cov}
\E\big[(X^\alpha)_{ab}\,(X^\beta)_{cd}\big]
=\frac{w(\alpha,\beta)}{N}\,\delta_{a d}\,\delta_{b c}.
\end{equation}

\begin{proof}[Proof of Proposition~\ref{prop:semi-closed-cov}]
Let $\mathcal P_2(p+q)$ be the set of pair partitions on $[p+q]$.
By Wick's formula,
\[
\E[W_1 W_2]
=\sum_{\pi\in\mathcal P_2(p+q)}\ \sum_{\mathbf i}\ 
\prod_{\{x,y\}\in \pi}\ \E\left[\left(X^{\tau(x)}\big)_{i_xi_{\rho(x)}}\ \big(X^{\tau(y)}\right)_{i_yi_{\rho(y)}}\right],
\]
where $\mathbf i=(i_1,\ldots,i_{p+q})$.
By \eqref{eq:elliptic-entry-cov}, we obtain the pairwise covariance for every pair $\{x,y\}\in\pi$:
\[
\E\left[\left(X^{\tau(x)}\right)_{i_x i_{\rho(x)}}\ \left(X^{\tau(y)}\right)_{i_y i_{\rho(y)}}\right]
=\frac{w \left(\tau(x), \tau(y)\right)}{N}\ \delta_{i_x i_{\rho(y)}}\,\delta_{i_{\rho(x)} i_y}.
\]
Hence,
\begin{equation}
\E[W_1 W_2]
=\sum_{\pi\in\mathcal P_2(p+q)}\ \sum_{\mathbf i}\ 
\prod_{\{x,y\}\in\pi}\left(\tfrac{w(\tau(x),\tau(y))}{N}\ \delta_{\,i_x i_{\rho(y)}}\delta_{i_{\rho(x)} i_y}\right).
\end{equation}

Partition $\mathcal P_2(p+q)$ into
\[
\mathcal P_2^{\times}(p+q)\coloneqq\Big\{\pi\in\mathcal P_2(p+q):\ \text{all pairs lie within the same cycle}\Big\},
\]
and its complement
\[
\mathcal P_2^{\mathrm{conn}}(p+q)\coloneqq\mathcal P_2(p+q)\setminus \mathcal P_2^{\times}(p+q)
=\Big\{\pi:\ \exists\,\{x,y\}\in\pi\ \text{with }x\le p<y\ \text{or}\ y\le p<x\Big\}.
\]
Note that the product 
\[\prod_{\{x,y\}\in\pi}w(\tau(x),\tau(y))\ \delta_{i_x i_{\rho(y)}}\delta_{i_{\rho(x)} i_y}\] 
can factorize into a product of a function of $(i_t)_{t=1}^p$ and a function of $(i_t)_{t=p+1}^{p+q}$, and the sum over indices splits accordingly. Via the natural bijection 
$\mathcal P_2^{\times}(p+q)\cong \mathcal P_2(p)\times\mathcal P_2(q)$, 
the two resulting sums are exactly the Wick expansions of $W_1$ and $W_2$. Therefore,
\[
\mathbb{E}[W_1] \mathbb{E}[W_2]
=\sum_{\pi\in\mathcal P_2^{\times}(p+q)}\ \sum_{\mathbf i}\ 
\prod_{\{x,y\}\in\pi}\Big(\tfrac{w(\tau(x),\tau(y))}{N}\ \delta_{i_x i_{\rho(y)}}\delta_{i_{\rho(x)} i_y}\Big),
\]
and the covariance equals the Wick contribution of the \emph{connected} pair partitions:
\begin{align}\label{eq:cov-linking}\notag
\kappa_2(W_1,W_2)&=\mathbb{E}[W_1W_2]-\mathbb{E}[W_1] \mathbb{E}[W_2]\\ \notag
&=\sum_{\pi\in\mathcal P_2^{\mathrm{conn}}(p+q)}\ \sum_{\mathbf i}\ 
\prod_{\{x,y\}\in\pi}\Big(\tfrac{w(\tau(x),\tau(y))}{N}\ \delta_{i_x i_{\rho(y)}}\delta_{i_{\rho(x)} i_y}\Big).
\end{align}

Rewriting the delta function for a pair $\{x,y\}$ with $\pi(x)=y$ as
\[
\delta_{\,i_x i_{\rho(y)}}\,\delta_{\,i_{\rho(x)}\, i_y}
=\delta_{\,i_x i_{\rho\pi(x)}}\,\delta_{\,i_{\rho\pi(y)} i_y}
\]
and multiplying over all pairs yields
\[
\prod_{\{x,y\}\in\pi}\delta_{\,i_x i_{\rho(y)}}\delta_{\,i_{\rho(x)} i_y}
=\prod_{x=1}^{p+q}\delta_{\,i_x i_{\rho\pi(x)}}.
\]
Thus, the delta functions force the indices to be constant along the cycles of $\rho\pi$; each cycle carries exactly one free index. In particular, the number of independent index sums equals the number of cycles, and hence
\[
\sum_{\mathbf i}\prod_{x=1}^{p+q}\delta_{\,i_x i_{\rho\pi(x)}}=N^{\,\#(\rho\pi)},
\quad
\text{where }\ \sum_{\mathbf i}=\sum_{i_1,\dots,i_{p+q}=1}^{N}.
\]
Therefore, after summing over indices, the contribution of $\pi$ is
\[
\Big(\tfrac{1}{N}\Big)^{(p+q)/2}N^{\,\#(\rho\pi)}
\prod_{\{x,y\}\in\pi} w\big(\tau(x),\tau(y)\big)
= N^{\,\#(\rho\pi)-(p+q)/2}\gamma^{\,s(\pi;\tau)}.
\]

Since $\pi$ links the two circles, the subgroup $\langle \pi,\rho\rangle$ acts transitively on $[p+q]$. 
By the Euler characteristic identity \cite[Theorem~5.9]{mingo2017free},
\[
\#(\pi)+\#(\pi^{-1}\rho)+\#(\rho)=(p+q)+2(1-g),
\]
where $g$ is the genus of $\pi$ relative to $\rho$ (the genus of $\pi$ relative to $\rho$ is the smallest $g$ such that the cycles of $\pi$ can be drawn on a surface of genus $g$). 
Here $\#(\pi)=(p+q)/2$ and $\#(\rho)=2$,
so
\[
\#(\rho\pi)=\#(\pi^{-1}\rho)=\frac{p+q}{2}-2g, \quad g\in\mathbb N_0.
\]
Consequently, \[N^{\,\#(\rho\pi)-(p+q)/2}=N^{-2g},\] and in the large-$N$ limit only the \emph{planar} ($g=0$) and connected pair partitions contribute. For such $\pi$, the Euler characteristic identity simplifies to
\[
\#(\pi)+\#(\pi^{-1}\rho)=p+q,
\]
and, together with connectedness, this is exactly the algebraic condition of the non-crossing annular pair partitions (see Remark~\ref{rmk:eqdef-ncannular}). Hence the contributing pair partitions are precisely $NC_2(p,q)$.
\end{proof}

\subsubsection{Real case}

The real case is quite different from the complex case. For  any $\alpha,\beta\in\{1,*\}$ (note that $X^{*}\coloneqq X^{T}$ for real $X$), we need to distinguish two types of Wick contractions, which we call cross and straight contractions. We encode them by the weights
\begin{equation}
\begin{split}
w_c(\alpha,\beta)=\gamma, w_s(\alpha,\beta)=1, & \; \text{if }\ \alpha=\beta,\\
w_c(\alpha,\beta)=1, w_s(\alpha,\beta)=\gamma, & \; \text{if }\ \alpha\neq\beta,
\end{split}
\end{equation}
where the subscripts c and s stand for cross and straight, respectively.

By Definition \ref{def:real-elliptic-matrix}, the entrywise covariance of $X$ is given by
\begin{equation}\label{eq:elliptic-entry-cov-real}
\E\left[(X^\alpha)_{ab}\,(X^\beta)_{cd}\right]
=\frac{1}{N}\left(
w_c(\alpha,\beta)\,\delta_{a d}\,\delta_{b c}
\ +\ 
w_s(\alpha,\beta)\,\delta_{a c}\,\delta_{b d}
\right).
\end{equation}

\begin{proposition}\label{prop:semi-closed-cov-real}
Let $X$ be the normalized real Gaussian elliptic matrix given in \eqref{eq:normalized-elliptic}, and let $W_1, W_2$ be given in \eqref{eq:W}. Then
\begin{equation}\label{eq:lim-cov-nc2-brief-real}
\lim_{N\to\infty}\kappa_2(W_1,W_2)
=\sum_{\pi\in NC_2(p,q)}
\left(
\gamma^{\,s(\pi;\tau)}
+\gamma^{s(\pi;\tau)+a(\pi)-2s_{sp}(\pi;\tau)}
\right),
\end{equation}
where $s(\pi;\tau):=\#\{\{x,y\}\in\pi:\ \tau(x)=\tau(y)\}$,  
$a(\pi)$ is the number of spokes of $\pi$, and 
$s_{\mathrm{sp}}(\pi;\tau)
:=\#\{\{x,y\}\in\pi:\ |\{x,y\}\cap[p]|=1,\ \tau(x)=\tau(y)\}$ is the number of same-type \emph{spokes} in $\pi$.
\end{proposition}

\begin{proof}
By \eqref{eq:elliptic-entry-cov-real}, we obtain
\begin{align}\notag
\mathbb{E}[W_1W_2]
=\sum_{\pi\in\mathcal P_2(p+q)} \sum_{\mathbf i}
\prod_{\{x,y\}\in\pi}\frac{1}{N}\Big[
&w_c(\tau(x),\tau(y))\delta_{i_x i_{\rho(y)}}\delta_{i_{\rho(x)} i_y}\\ \notag
&\quad+w_s(\tau(x),\tau(y))\delta_{i_x i_y}\delta_{i_{\rho(x)} i_{\rho(y)}}
\Big],
\end{align}
where $\mathbf i=(i_1,\ldots,i_{p+q})$.

Similar to the complex case, pair partition $\pi$ whose pairs lie within a single circle 
produce a factorized index sum and reproduce $\mathbb{E}[W_1]\mathbb{E}[W_2]$.
Hence
\begin{equation}\label{eq:real-assignments}
\kappa_2(W_1,W_2)
=\sum_{\pi\in\mathcal P_2^{\mathrm{conn}}(p+q)}\ \sum_{\mathbf i}\
\prod_{\{x,y\}\in\pi}\frac{1}{N}\Big(
w_c\,\delta_{i_x i_{\rho(y)}}\delta_{i_{\rho(x)} i_y}
+w_s\,\delta_{i_x i_y}\delta_{i_{\rho(x)} i_{\rho(y)}}
\Big).
\end{equation}

\noindent\emph{Step 1: Power counting and the leading order.}
Denote the cross and straight constraints by
\[
\delta^{c}_{x,y}\coloneqq \delta_{i_x i_{\rho(y)}}\delta_{i_{\rho(x)} i_y},
\quad
\delta^{s}_{x,y}\coloneqq \delta_{i_x i_y}\delta_{i_{\rho(x)} i_{\rho(y)}},
\]
such that for each $\{x,y\}\in\pi$ the bracket in \eqref{eq:real-assignments} reads $w_c\,\delta^{c}_{x,y}+w_s\,\delta^{s}_{x,y}$.
Expanding the product over $\pi$ amounts to summing over assignments 
$\varepsilon:\pi\to\{c,s\}$:
\[
\prod_{\{x,y\}\in\pi}\Big(w_c\,\delta^{c}_{x,y}+w_s\,\delta^{s}_{x,y}\Big)
=
\sum_{\varepsilon:\,\pi\to\{c,s\}}
 \prod_{\{x,y\}\in\pi} w_{\varepsilon(\{x,y\})}
 \prod_{\{x,y\}\in\pi}\delta^{\,\varepsilon(\{x,y\})}_{x,y}.
\]
For a fixed assignment $\varepsilon$, we associate an index permutation $\theta_\varepsilon\in S_{p+q}$ such that 
\begin{equation}\label{def:index-permut}
\prod_{\{x,y\}\in\pi}\delta^{\,\varepsilon(\{x,y\})}_{x,y}
=
\prod_{t=1}^{p+q}\delta_{\,i_t i_{\theta_\varepsilon(t)}}.
\end{equation}
Namely, $\theta_\varepsilon$ is an index permutation if and only if the delta functions imply $i_t=i_{\theta(t)}$ for every $t\in[p+q]$; that is, the indices are constant along the cycles of $\theta$. The choice of $\theta_\varepsilon$ is not unique (the cyclic order in each cycle may change), but its cycle count $\#(\theta_\varepsilon)$ is uniquely determined by the assignment and is the only quantity relevant for index-sum power counting. Therefore,
\[
\sum_{\mathbf i}\prod_{t=1}^{p+q}\delta_{\,i_t i_{\theta_\varepsilon(t)}}
= N^{\,\#(\theta_\varepsilon)}.
\]
Then the contribution of a fixed
$\pi$ and assignment $\varepsilon$ is
\[
N^{\,\#(\theta_\varepsilon)-(p+q)/2}\ \prod_{\{x,y\}\in\pi} w_{\varepsilon(\{x,y\})}.
\]

We now discuss different assignments:
\begin{enumerate}
	\item \emph{All-cross assignment.} If every pair $\{x,y\}\in\pi$ is realized by the cross
constraint $\delta^{c}_{x,y}$ (i.e.\ $\varepsilon\equiv c$),
this coincides exactly with the complex case; in particular, $\theta_\varepsilon=\rho\pi$.
In this case, only the non-crossing annular pair partitions $\pi\in NC_2(p,q)$ contribute after taking the limit.
    \item \emph{Any straight arc pair is subleading.}
Fix an arc pair $\{x,y\}$ on one circle and assume that $y$ is encountered after $x$ along $\rho$.
Let $I_1\coloneqq(x,y)_\rho$ and $I_2\coloneqq(y,x)_\rho$ be the open $\rho$-intervals on that circle.
By non-crossing, these intervals are $\pi$-invariant: $\pi(I_1)=I_1$ and $\pi(I_2)=I_2$.
Set $S_1:=I_1\cup\{y\}$ and $S_2:=I_2\cup\{x\}$. Then for any $a\in I_2$ we have
$\rho\pi(a)\in S_2$, and also $\rho\pi(x)=\rho(y)\in I_2$, hence $\rho\pi(S_2)\subseteq S_2$.
By induction, $(\rho\pi)^k(x)\in S_2$ for all $k\ge0$. Since $y\notin S_2$, it follows that
$(\rho\pi)^k(x)\neq y$ for all $k\ge1$. Therefore $x$ and $y$ lie in \emph{distinct} cycles of
$\rho\pi$. Thus, realizing the straight constraint $\delta^{s}_{x,y}$ on this arc will identify two previously independent index cycles and \emph{decreases} the cycle count: $\#(\theta_\varepsilon)=\#(\rho\pi)-1$. Thus any assignment with at least one arc taken straight is suppressed by a factor $N^{-1}$, and is subleading in the large-$N$ limit. More generally, if $m$ arc pairs with pairwise disjoint endpoints are taken straight, then $\#(\theta_\varepsilon)=\#(\rho\pi)-m$ and the contribution is suppressed by $N^{-m}$.
\item \emph{Mixed assignment: Spokes straight, arcs pairs cross.}
On the outer circle perform the following global bijective relabeling of dummy indices:
\[
j_y\coloneqq i_{\rho(y)}\quad\text{for all }y\in p+[q],
\]
Since $y\mapsto\rho(y)$ is a permutation on the outer circle, this relabeling is a bijection and
does not change the total index sum. Now reverse the cyclic order of the trace permutation on the outer circle while keeping the order on the inner circle unchanged, namely define
\begin{equation}\label{eq:hat-rho}
\hat\rho(z)\coloneqq
\begin{cases}
\rho(z),& z\in [p]\ \text{(inner circle)},\\[2pt]
\rho^{-1}(z),& z\in p+[q]\ \text{(outer circle)}.
\end{cases}
\end{equation}
For a spoke $\{x,y\}$ with $y$ on the outer circle, the straight constraints \(i_x=i_y,\ i_{\rho(x)}=i_{\rho(y)}\) become, under the relabeling,
\[
i_x=j_{\rho^{-1}(y)}=j_{\hat\rho(y)},\qquad i_{\hat \rho(x)}=j_y,
\]
which are exactly the cross constraints with respect to $\hat\rho$. For an outer-circle arc pair $\{u,v\}\subset p+[q]$ the cross constraints read \(i_u=i_{\rho(v)},\ i_{\rho(u)}=i_v.\) Substituting $j_y=i_{\rho(y)}$ yields
\[
j_{\rho^{-1}(u)}=j_v,\qquad j_u=j_{\rho^{-1}(v)},
\]
which is equivalent to
\(
j_u=j_{\hat\rho(v)},\ j_{\hat\rho(u)}=j_v,
\)
i.e. still cross with respect to $\hat\rho$. Indices and the cyclic order of trace permutation remain unchanged on the inner circle, so inner-circle arc pairs stay cross under this relabeling.

With this replacement, every straight spoke turns into a cross spoke and every cross arc pairs remains cross (relative to $\hat\rho$). Hence the mixed assignment coincides with the complex case (with $\rho$ replaced by $\hat\rho$), so it contributes at the leading order; its contributions are exactly those $\pi\in NC_2(p,q)$ with the outer circle read in reversed cyclic order. 

\item \emph{Mixed spoke assignments are subleading.}
At a spoke $\{x,y\}$, the cross constraints carry two distinct index variables $i_x$ and $i_y$ (one propagated through $\delta_{i_x i_{\rho(y)}}$, the other through $\delta_{i_{\rho(x)} i_y}$). If this single spoke is switched to the straight constraint, the equalities $i_x=i_y$ and $i_{\rho(x)}=i_{\rho(y)}$ identify those two variables and reduce the number of independent index cycles by one; that is, $\#(\theta_\varepsilon)= \#(\rho\pi)-1= \#(\hat \rho \pi)-1$ relative to the all spokes cross or all spokes straight assignments. Hence any assignment in which some spokes are cross and some are straight (with arcs kept cross) does not contribute at leading order.
\end{enumerate}
Above all, in the large-$N$ limit exactly two assignments contribute at leading order:
\begin{itemize}
\item \underline{\emph{All pairs cross.}} The index permutation $\theta_\varepsilon=\rho\pi$, so the contributing diagrams are precisely the standard non-crossing annular pair partitions $\pi\in NC_2(p,q)$.
\item \underline{\emph{Spokes straight, arc pairs cross.}} The index permutation $\theta_\varepsilon=\hat\rho\pi$, which again yields non-crossing annular pair partitions---this is the \emph{outer-reversed} counterpart of the standard non-crossing annular pair partitions.
\end{itemize}

\noindent
\emph{Step 2: Weights of the two assignments.}

\noindent
\begin{enumerate}

\item \emph{All pairs cross.} Each same-type pair (whether arc or spoke) contributes a factor $\gamma$, each mixed-type pair contributes $1$. Therefore the total weight is
\[
\prod_{\{x,y\}\in\pi} w_c(\tau(x),\tau(y))=\gamma^{\,s(\pi;\tau)}.
\]

\item \emph{Spokes straight, arc pairs cross.} On spokes, the straight constraint swaps the $1/\ast$ roles: a same-type spoke now contributes weight $1$ while a different-type spoke contributes $\gamma$. Hence the spoke factor equals $\gamma^{\,a(\pi)-s_{sp}(\pi;\tau)}$. For arc pairs the contribution is the same as in the complex case and  contributes $\gamma$ exactly for the same-type arc pairs. However, the global type word on the outer circle should be read in reversed cyclic order. Together with the spoke factor this yields the total
\[
\gamma^{\,s(\pi;\tau)-s_{sp}(\pi;\tau)}\cdot
\gamma^{\,a(\pi)-s_{sp}(\pi;\tau)}
=\gamma^{\,s(\pi;\tau)+a(\pi)-2s_{sp}(\pi;\tau)}.
\]
\end{enumerate}
Collecting the preceding cases into \eqref{eq:real-assignments} and summing the two leading contributions over $\pi\in NC_2(p,q)$ yields \eqref{eq:lim-cov-nc2-brief-real}.
\end{proof}

\begin{remark}\label{spoke-straight=outer-reversed-transposed}
Recall that we work in the two-circle notation: the inner cycle is $[p]$, the outer cycle is $p+[q]$,
and the trace permutation is
\[
  \rho = (1\,2\,\cdots\,p)\,(p+1\,\cdots\,p+q).
\]
Accordingly, for any two trace products of matrices $W_{\mathrm{in}}=\Tr(A_1\cdots A_p)$ and
$W_{\mathrm{out}}=\Tr(A_{p+1}\cdots A_{p+q})$ we encode the product of traces as
\[
  \Tr(A_1\cdots A_p)\cdot \Tr(A_{p+1}\cdots A_{p+q})
  =
  \sum_{i_1,\ldots,i_{p+q}}
  \ \prod_{t=1}^{p+q} (A_t)_{\,i_t\,i_{\rho(t)}},
\]
so that the successor along each circle is precisely the index contraction inside the corresponding trace.

Consider the mixed assignment in which all spokes are straight while all arc pairs
are cross. We claim that, at the level of index contractions, the mixed assignment for the \emph{real} ensemble
coincides with the \emph{complex} all-cross channel after replacing the outer product by its transpose
read in reversed order; that is, with $(W_{\mathrm{in}},W_{\mathrm{out}}^{T})$ in the complex ensemble.
Precisely, define
\[
  W_{\mathrm{out}}^{T}
  \coloneqq
  \Tr\big((A_{p+1}\cdots A_{p+q})^T\big)
  = \Tr(A_{p+q}^T\cdots A_{p+1}^T).
\]
Then the mixed contribution equals the complex all-cross channel for $(W_{\mathrm{in}},W_{\mathrm{out}}^{\,T})$:
\[
  \sum_{\pi\in NC_2(p,q)} \gamma^{\,s(\pi;\tau)+a(\pi)-2\,s_{\mathrm{sp}}(\pi;\tau)}
  =
  \lim_{N\to\infty}\kappa_2^{\mathbb C}\big(W_{\mathrm{in}},\,W_{\mathrm{out}}^{\,T}\big),
\]
while the all-cross contribution equals $\lim_{N\to\infty}\kappa_2^{\mathbb C}(W_{\mathrm{in}},W_{\mathrm{out}})$.
Adding the two channels yields the relation
\[
   \lim_{N\to\infty}\kappa_2^{\mathbb R}\big(W_{\mathrm{in}},W_{\mathrm{out}}\big)
  =
  \lim_{N\to\infty}\kappa_2^{\mathbb C}\big(W_{\mathrm{in}},W_{\mathrm{out}}\big)
  +
  \lim_{N\to\infty}\kappa_2^{\mathbb C}\big(W_{\mathrm{in}},W_{\mathrm{out}}^{T}\big),
\]
where the superscripts $\mathbb R$ and $\mathbb C$ indicate the underlying ensemble (real vs.\ complex elliptic matrices).

Indeed, expand $W_{\mathrm{out}}$ as
\[
  \Tr(A_{p+1}\cdots A_{p+q})
  =
  \sum_{j_{p+1},\ldots,j_{p+q}}
  \ \prod_{t=p+1}^{p+q} (A_t)_{\,j_t,\,j_{\rho(t)}}.
\]
Recall that the mixed assignment coincides with the complex case with $\rho$ replaced by $\hat \rho$. We make the bijective relabeling on the outer circle $i_t:=j_{\rho(t)}$ for $t\in p+[q]$.
Then for each $t\in p+[q]$,
\[
  (A_t)_{\,j_t,\,j_{\rho(t)}}
  =
  (A_t)_{\,i_{\rho^{-1}(t)},\,i_t}
  =
  (A_t^T)_{\,i_t,\,i_{\hat\rho(t)}},
\]
where $\hat\rho$ is defined in \eqref{eq:hat-rho}.
Summing over the relabeled indices yields
\[
  \sum_{i_{p+1},\ldots,i_{p+q}}
  \prod_{t=p+1}^{p+q} (A_t^T)_{\,i_t,\,i_{\hat\rho(t)}}
  = \Tr(A_{p+q}^T\cdots A_{p+1}^T).
\]

This matches the phenomenon described in \cite{Redelmeier2014}: annular spoke diagrams with two oppositely oriented circles of the annulus, in which the matrix transpose appears.
\end{remark}

\subsection{Limiting covariance via arc weights}
Fix a type word $\tau: [p+q] \to \{1, \ast\}$ and $\pi\in NC_2(p,q)$. In the spoke-arc decomposition, pair weights $w(\tau(x),\tau(y))\in\{1,\gamma\}$ factor multiplicatively over spokes and arcs. The arc contributions will be packaged into local \emph{arc weights} $F(\cdot)$ defined below; thus, at the spoke-arc level, the only remaining $\gamma$-dependence comes from how many spokes are same-type. This motivates the following definition for a labeled spoke-arc configuration.

\begin{definition}[Number of same-type spokes]\label{def:No-of-sametype-spokes}
Given a type word $\tau$ and a labeled spoke-arc configuration $(a,\mathbf o,\boldsymbol\iota,U,V;\ \{\pi^{\mathbf{o}}_r\}_{r=1}^a,\ \{\pi^{\boldsymbol\iota}_r\}_{r=1}^a)$ in the spoke-arc decomposition, define its number of same-type spokes by
\begin{equation}
s_{\mathrm{sp}}(a;\tau;U,V) \coloneqq \#\big\{\,r\in\{1,\dots,a\}:\ \tau(u_r)=\tau(v_r)\,\big\}.
\end{equation}
\end{definition}

\begin{remark}
This quantity is well-defined on the \(\mathbb Z_a\)-orbit of \((U,V)\), i.e., \(s_{\mathrm{sp}}(a;\tau;U,V)\allowbreak=\allowbreak s_{\mathrm{sp}}(a;\tau;U^{(s)},V^{(s)})\) for all simultaneous cyclic relabelings $s\in \{0,\ldots,a-1\}$, where $U^{(s)}=(u_{1+s},\ldots,u_{a+s})$ and $V^{(s)}=(v_{1+s},\ldots,v_{a+s})$.
\end{remark}

For fixed $(a,\mathbf o,\boldsymbol\iota)$, the spoke endpoint lists $U\subset[p]$ and $V\subset p+[q]$ are determined by their starting points: a choice of start on $[p]$ (resp.\ on $p+[q]$) produces exactly $p$ (resp.\ $q$) labelings realizing the same arc lengths; taking the quotient by the simultaneous cyclic relabeling of the $a$ spokes removes the overcounting. This yields the following lemma, which counts the labeling of spokes. 

\begin{lemma}\label{lem:abs-spoke-count}
Fix admissible data \((a,\mathbf o,\boldsymbol\iota)\),
and let \(\mathcal U_{a,\mathbf o}\) (resp.\ \(\mathcal V_{a,\boldsymbol\iota}\))
denote the set of ordered outer (resp.\ inner) endpoint lists
\(U=(u_1,\dots,u_a)\subset [p]\) (resp.\ \(V=(v_1,\dots,v_a)\subset p+[q]\))
written in cyclic order and realizing the arc lengths \(\mathbf o\) (resp.\ \(\boldsymbol\iota\)).
Then
\[
|\mathcal U_{a,\mathbf o}|=p,\quad |\mathcal V_{a,\boldsymbol\iota}|=q.
\]
Moreover, modulo the cyclic relabeling by \(\mathbb Z_a\), the number of distinct spoke labelings equals
\[
\frac{|\mathcal U_{a,\mathbf o}|\,|\mathcal V_{a,\boldsymbol\iota}|}{|\mathbb Z_a|}
=\frac{pq}{a}.
\]
\end{lemma}

\begin{proof}
Fix \(\mathbf o\). Given any choice of \(u_1\in[p]\), define \(U\) recursively by
\[
u_{r+1}= u_r+o_r+1 \pmod{p},\quad (r=1,\dots,a;\ \text{indices modulo }a),
\]
where ``\(+1\)'' corresponds to the next spoke endpoint after traversing an outer arc of length \(o_r\).
Because \(\sum_{r=1}^a(o_r+1)=p\), this construction yields an ordered \(a\)-tuple of distinct
positions in \([p]\) and closes after \(a\) steps. Conversely, any \(U\) realizing \(\mathbf o\)
arises uniquely from its first entry \(u_1\). Hence \(|\mathcal U_{a,\mathbf o}|=p\).
The same argument on the inner circle with \(\boldsymbol\iota\) gives
\(|\mathcal V_{a,\boldsymbol\iota}|=q\).

Two pairs \((U,V)\) and \((U',V')\) yield the same spoke labeling if and only if they differ by a
common cyclic relabeling \(r\mapsto r+s\) for some \(s\in\{0,\dots,a-1\}\). Therefore the number of distinct spoke labelings is $pq /a.$
\end{proof}

By Lemma~\ref{lem:abs-spoke-count}, the only remaining freedom in a labeled spoke-arc configuration is the independent choice of non-crossing pair partitions within each arc. It is convenient to package the total contribution of a single arc into a local quantity that depends only on the type word carried by that arc (and not on its absolute position or on the global spoke labeling). This motivates the following definition.

\begin{definition}[Arc weight]\label{def:arc-weight-F}
Let $\tau: [2n] \to \{1,\ast\}$ be an even-length type word placed along a single arc (with $2n$ marked points in cyclic order). Define the \emph{arc weight}
\[
F(\tau)\coloneqq\sum_{\pi\in NC_2(2n)}\gamma^{\,s(\pi;\tau)},
\]
where $\gamma\in[-1,1]$ is the elliptic parameter.
\end{definition}

\begin{example}\label{rmk:mc-arc-weight}
Let $\tau: [2n] \to \{1, \ast\},$ and let $C_n$ be the $n$-th Catalan number.
\begin{itemize}
\item If $\tau$ is pure, then $F(\tau)=\gamma^{n}C_n;$ If $\tau$ is alternating, then $F(\tau)=C_n.$
\item Let $X$ be the normalized elliptic Gaussian matrix. By Theorem \ref{thm:moment} (see Adhikari--Bose \cite[Lemma~1]{adhikari2019brown}), the arc weight is exactly the limiting first moment on that arc:
\begin{equation}\label{id-arcweight-moment}
F(\tau)=\lim_{N\to\infty}\mathbb E\tr\left(X^{\tau(1)}\cdots X^{\tau(2n)}\right).
\end{equation}
\end{itemize}
\end{example}

Fix a labeled spoke-arc configuration and let $\tau$ be the global type word on $[p+q]$. The underlying pair partition $\pi$ is the disjoint union of (i) the $a$ spokes and (ii) the arc pairs on each outer and inner arc. Since the same-type indicator depends only on the two endpoints of a pair, the total same-type count splits additively as
\[
s(\pi;\tau)
=s_{\mathrm{sp}}(\pi;\tau)
+\sum_{r=1}^{a} s\big(\pi^{\mathbf o}_r; \tau^{\mathbf o}_r\big)
+\sum_{r=1}^{a} s\big(\pi^{\boldsymbol \iota}_r; \tau^{\boldsymbol \iota}_r\big),
\]
where $\tau^{\mathbf o}_r$ (resp.\ $\tau^{\boldsymbol\iota}_r$) denotes the restriction of $\tau$ to the
$r$-th outer (resp.\ inner) arc.
Consequently the $\gamma$-weight separates multiplicatively into a spoke factor and independent arc factors,
which yields the following simple lemma.

\begin{lemma}[Weight factorization]\label{lem:weightfac}
Let $\pi\in NC_2(p,q)$ be represented by a labeled spoke-arc configuration
\(
(a,\mathbf o,\boldsymbol\iota, U,V;\ \{\pi^{\mathbf o}_r\}_{r=1}^a,\ \{\pi^{\boldsymbol\iota}_r\}_{r=1}^a)
\).
Then, with respect to the global type word $\tau$, we have 
\begin{equation}
\gamma^{\,s(\pi;\tau)} =\;\gamma^{\,s_{\mathrm{sp}}(\pi;\tau)}
   \prod_{r=1}^{a} \gamma^{\,s(\pi^{\mathbf o}_r; \tau^{\mathbf o}_r)}
   \prod_{r=1}^{a} \gamma^{\,s(\pi^{\boldsymbol \iota}_r; \tau^{\boldsymbol \iota}_r)}.
\end{equation}
\end{lemma}

Combining the above ingredients, and introducing the local arc contributions via the single-arc weights $F(\cdot)$, then summing over all admissible data $(a,\mathbf o,\boldsymbol\iota)$ and over the $\mathbb Z_a$-orbits of endpoint lists $(U,V)$, we obtain the following limiting covariance formula for the complex case.

\begin{corollary}\label{cor:spoke-arc}
Given any $p,q\ge1$ and type word $\tau: [p+q] \to \{1, \ast\}$. Let $X$ be the normalized complex Gaussian elliptic matrix given in \eqref{eq:normalized-elliptic}, and let $W_1, W_2$ be given in \eqref{eq:W}. We have
\begin{equation}\label{eq:spoke-arc-complex}
\lim_{N\to\infty}\kappa_2(W_1,W_2)
=\sum_{a=1}^{\min\{p,q\}}
 \sum_{\mathbf o,\boldsymbol{\iota}}
 \sum_{[(U,V)]\in \mathcal O_{a,\mathbf o,\boldsymbol\iota}}
\gamma^{\,s_{\mathrm{sp}}(a;\tau;U,V)}
\prod_{r=1}^{a} F\big(\tau^{\mathbf o}_r\big)\ \prod_{r=1}^{a} F\big(\tau^{\boldsymbol\iota}_r\big),
\end{equation}
where \(\mathcal O_{a,\mathbf o,\boldsymbol\iota}\) denotes the set of \(\mathbb Z_a\)-orbits of ordered
endpoint lists \((U,V)\) realizing the data \((a,\mathbf o,\boldsymbol\iota)\). 
In particular, \(|\mathcal O_{a,\mathbf o,\boldsymbol\iota}|=\frac{pq}{a}\). 
\end{corollary}

\begin{proof}
Recall that every $\pi\in NC_2(p,q)$ can be encoded by a labeled spoke-arc configuration
\[
(a,\mathbf o,\boldsymbol\iota,U,V,\{\pi_r^{\mathbf o}\}_{r=1}^a,\{\pi_r^{\boldsymbol\iota}\}_{r=1}^a).
\]
By Lemma~\ref{lem:weightfac}, the pair parition $\pi$ splits as the disjoint union
\[
\pi = \{\{u_r,v_r\}\}_{r=1}^{a}
\ \sqcup\ \bigsqcup_{r=1}^{a}\pi^{\mathbf o}_r
\ \sqcup\ \bigsqcup_{r=1}^{a}\pi^{\boldsymbol\iota}_r,
\] 
and the total $\gamma$-weight of the fixed labeled spoke-arc configuration equals
\[
\gamma^{\,s_{\mathrm{sp}}(a;\tau;U,V)}
\prod_{r=1}^{a} F\big(\tau^{\mathbf o}_r\big)\;
\prod_{r=1}^{a} F\big(\tau^{\boldsymbol\iota}_r\big).
\]

Summing the above contribution over all admissible $(a,\mathbf o,\boldsymbol\iota)$ and over all $[(U,V)]\in\mathcal O_{a,\mathbf o,\boldsymbol\iota}$ yields Equation \eqref{eq:spoke-arc-complex}. Finally, by Lemma~\ref{lem:abs-spoke-count}, for each admissible $(a,\mathbf o,\boldsymbol\iota)$ there are exactly $|\mathcal O_{a,\mathbf o,\boldsymbol\iota}|=\frac{pq}{a}$ distinct $\mathbb Z_a$-orbits of endpoint lists $(U,V)$.
\end{proof}

\begin{corollary}[The real case]\label{cor:real-spoke-arc}
Given any $p,q\ge1$ and type word $\tau: [p+q] \to \{1, \ast\}$. Let $X$ be the normalized real Gaussian elliptic matrix given in \eqref{eq:normalized-elliptic}, and let $W_1, W_2$ be given in \eqref{eq:W}. We have
\begin{equation}\label{eq:spoke-arc-real}
\begin{aligned}
		&\lim_{N\to\infty}\kappa_2(W_1,W_2)\\
&=\left[(\text{complex})+\sum_{a=1}^{\min\{p,q\}}
 \sum_{\mathbf o,\boldsymbol{\iota}}
 \sum_{[(U,V)]\in \mathcal O_{a,\mathbf o,\boldsymbol\iota}}
\gamma^{a-s_{sp}(a;\tau;U,V)}
\prod_{r=1}^{a}F\left( \tau^{\mathbf o}_r \right)\prod_{r=1}^{a}F\left(\tau^{\boldsymbol\iota}_r\right)\right].
\end{aligned}
\end{equation}
Here ``(complex)'' denotes the right-hand side of~\eqref{eq:spoke-arc-complex}.
\end{corollary}

\begin{proof}
This is an immediate combination of the ingredients already established.
By Proposition~\ref{prop:semi-closed-cov-real} and the spoke-arc decomposition, each
$\pi\in NC_2(p,q)$ corresponds to a $\mathbb Z_a$-orbit of labeled spoke-arc configurations
$(a,\mathbf o,\boldsymbol\iota,U,V,(\pi^{\mathbf o}_r),(\pi^{\boldsymbol\iota}_r))$.
By Lemma~\ref{lem:weightfac}, the contribution of such a configuration
is the sum of the two pieces:
\[
\underbrace{\gamma^{\,s_{sp}(a;\tau;U,V)}\prod_{r=1}^a F(\tau^{\mathbf o}_r)\prod_{r=1}^a F(\tau^{\boldsymbol\iota}_r)}_{\text{all pairs cross}}
\;+\;
\underbrace{\gamma^{\,a-s_{sp}(a;\tau;U,V)}\prod_{r=1}^a F\!\big((\tau^{\mathbf o}_r)^{\mathrm{rev}}\big)\prod_{r=1}^a F(\tau^{\boldsymbol\iota}_r)}_{\text{spokes straight,\ arcs cross (outer reversed)}}.
\]
Here $(\cdot)^{\mathrm{rev}}$ denotes reading the outer arc type word in the reversed cyclic order; in particular $F((\tau^{\mathbf o}_r)^{\mathrm{rev}})=F(\tau^{\mathbf o}_r)$, thus reversing the outer cyclic order does not affect the arc weight.
\end{proof}

\begin{remark}[Multi-color extension]\label{rmk:multi-color}
Suppose that the letters on the inner/outer type words carry colors in a finite set $\mathcal C$, encoded by
$c:\{1,\dots,|\tau|\}\to\mathcal C$ (for a single arc of length $2n$, this restricts to $c:\{1,\dots,2n\}\to\mathcal C$).
For a \emph{multi-color} arc, we restrict to color-respecting pair partitions
\begin{equation}
NC_2^{c}(2n)\coloneqq\Big\{\pi\in NC_2(2n):\ c(x)=c(y)\ \text{for all }\{x,y\}\in\pi\Big\},
\end{equation}
and define the \emph{multi-color arc weight}
\begin{equation}
F_{c}(\tau)\coloneqq\sum_{\pi\in NC_2^{c}(2n)}\gamma^{\,s(\pi;\tau)}.
\end{equation}
It reduces to the single-color weight when $c$ is constant: $F_c(\tau)=F(\tau)$. Moreover, by \cite[Lemma~9]{mingo2017free},
\[
F_{c}(\tau)=\lim_{N\to\infty}\E\tr\left(\prod_{t=1}^{2n}\left(X^{(c(t))}\right)^{\tau(t)}\right).
\]

The spoke-arc decomposition admits a multi-color refinement by restricting to color-respecting spoke labelings. Define
\begin{align}\notag
\mathcal O^{\,c}_{a,\mathbf o,\boldsymbol\iota}
\coloneqq\Big\{[(U,V)]\in\mathcal O_{a,\mathbf o,\boldsymbol\iota}\ \Big|\ \exists\,k\in\mathbb Z_a\ \text{s.t. } 
&\forall r\in\mathbb Z_a,\ c(u_r)=c(v_{k-r})\\ \notag
\text{or } \,&\forall r\in\mathbb Z_a,\ c(u_r)=c(v_{k+r})\Big\},
\end{align}
i.e., the subset of \(\mathcal O_{a,\mathbf o,\boldsymbol\iota}\) such that every spoke is color-respecting (its two endpoints carry the same color). Then the real-case formula \eqref{eq:spoke-arc-real} becomes
\begin{equation}\label{eq:mc-spoke-arc-real}
\begin{aligned}\notag
\lim_{N\to\infty}\kappa_2(W_1,W_2)
&=\text{(\emph{complex})}
+\sum_{a=1}^{\min\{p,q\}}\ \sum_{\mathbf o,\boldsymbol\iota}\ 
 \sum_{[(U,V)]\in \mathcal O^c_{a,\mathbf o,\boldsymbol\iota}}
\gamma^{\,a-s_{\mathrm{sp}}(a;\tau;U,V)}\\[-2pt]
&\hspace{3.3cm}\times
\prod_{r=1}^{a}F_c\big((\tau^{\mathbf o}_r)^{\mathrm{rev}}\big)\ \prod_{r=1}^{a}F_c\big(\tau^{\boldsymbol\iota}_r\big),
\end{aligned}
\end{equation}
where the \emph{ ``(complex)'' } term is obtained from \eqref{eq:spoke-arc-complex} by replacing each
$\mathcal O_{a,\mathbf o,\boldsymbol\iota}$ with $\mathcal O^{\,c}_{a,\mathbf o,\boldsymbol\iota}$ and
each $F(\cdot)$ with $F_c(\cdot)$. 
\end{remark}

\section{Proof of main results}\label{sec:proof}

\subsection{Proof of Theorem~\ref{thm:spoke-arc-closed} and Theorem~\ref{thm:spoke-arc-closed-real}}

Fix integers \(p,q\ge1\). For any $\pi \in NC_2(p, q)$, let the data $(a,\mathbf o,\boldsymbol\iota,U,V;\ \{\pi^{\mathbf{o}}_r\}_{r=1}^a,\ \{\pi^{\boldsymbol\iota}_r\}_{r=1}^a)$ be its spoke-arc decomposition given in Definition \ref{def:spoke-arc-decomposition}.
Suppose that $p-a$ and $q-a$ are even numbers. After removing the spokes, each circle is cut into $a$ open arcs. Denote the remaining
points on the $a$ outer (resp.\ inner) arcs by
\[
(o_1,\dots,o_a)\quad\text{and}\quad (\iota_1,\dots,\iota_a),
\]
respectively. Then it is clear that $o_r, \iota_r, r=1, \ldots, a$ are even numbers and 
\begin{equation*}
\sum_{r=1}^a o_r=p-a, \;\;  \sum_{r=1}^a \iota_r=q-a.
\end{equation*}

Let \(\tau: [p+q] \to \{1, \ast\}\) be a type word on \([p+q]\), and write 
$
\tau^{o}_r (\text{resp.\ }\tau^{\iota}_r)
$
for its restriction to the \(r\)-th outer (resp.\ inner) arc. Recall the \emph{arc weight}
\[
F(\tau)\ \coloneqq \sum_{\pi\in NC_2(|\tau|)} \gamma^{\,s(\pi;\tau)}\quad
\big(\,|\tau|\ \text{even}\,\big),
\]
where $|\tau|$ denotes the length of $\tau.$

\begin{lemma}\label{lem:arc-compress-unified}
With the above-mentioned notations and assumptions, one has:
\begin{align}
\sum_{\substack{o_1+\cdots+o_a=p-a\\ o_r\in 2\mathbb N_0}}
\ \prod_{r=1}^a F\big(\tau^{\mathbf o}_r\big)
&=
\begin{cases}
\displaystyle \gamma^{\frac{p-a}{2}}\ \mathrm{FC}\left(a, \frac{p-a}{2}\right),
& \tau^{\mathbf o}_r, r=1, \ldots, a, \; \text{are pure},\\[10pt]
\displaystyle \mathrm{FC}\left(a, \frac{p-a}{2}\right),
& \tau^{\mathbf o}_r, r=1, \ldots, a, \text{are alternating}.
\end{cases}
\label{eq:outer-compress-F}
\\[6pt]
\sum_{\substack{\iota_1+\cdots+\iota_a=q-a\\ \iota_r\in 2\mathbb N_0}}
\ \prod_{r=1}^a F\big(\tau^{\boldsymbol{\iota}}_r\big)
&=
\begin{cases}
\displaystyle \gamma^{\frac{q-a}{2}}\  \mathrm{FC}\left(a, \frac{q-a}{2}\right),
& \tau^{\boldsymbol{\iota}}_r, r=1, \ldots, a, \; \text{are pure},\\[10pt]
\displaystyle \mathrm{FC}\left(a, \frac{q-a}{2}\right),
& \tau^{\boldsymbol{\iota}}_r, r=1, \ldots, a, \text{are alternating},
\end{cases}
\label{eq:inner-compress-F}
\end{align}
where $\mathrm{FC}(a,n)\coloneqq \frac{a}{2n+a}\binom{2n+a}{n}$ is the Fuss-Catalan number.
\end{lemma}

\begin{proof}
We only prove \eqref{eq:outer-compress-F}, and the proof of \eqref{eq:inner-compress-F} is similar. 
For a formal power series $f(z)=\sum_{n\ge0} a_n z^n$ we write
\[
[z^n]f(z)\coloneqq a_n.
\]
We firstly claim the following simple result, which is a consequence of the Lagrange inversion theorem (see, e.g., \cite[Appendix~A.6]{FlajoletSedgewick2009} and the recent inductive proof in \cite{SuryaWarnke2023}). For integers $a\ge 1$ and $n\ge 0$, we have
\begin{equation}\label{eq:Lagrange-inversion}
[z^{n}]\,\mathcal C(z)^{\,a}
\;=\;\frac{a}{2n+a}\binom{2n+a}{n} = \mathrm{FC}(a,n),
\end{equation}
where $\mathcal C(z) = \frac{1-\sqrt{1-4z}}{2z}$ is the generating function of Catalan numbers and $\mathrm{FC}(a,n)$ is called the Fuss-Catalan number (or the Raney number \cite{Raney1960}) for $a, n \in \mathbb{N}$. Note that $\mathcal C(z)=1+z\,\mathcal C(z)^2$ (see, e.g., \cite[Lemma~2.21]{nica2006lectures}).
Set $y(z):=\mathcal C(z)$ and $x(z):=y(z)-1$. Then $x(z)$ solves the Lagrange form
\(
x(z)=z\,\Phi\big(x(z)\big)\) with \(\Phi(u)=(u+1)^2.
\)
Since \(y(z)=x(z)+1\), we have $[z^n]\,y(z)^{a}=[z^n]\,(x(z)+1)^{a}$. By the Lagrange inversion theorem,
for $n\ge 1$ and any differentiable $g$,
\[
[z^{n}]\,g\big(x(z)\big)\;=\;\frac{1}{n}\,[u^{\,n-1}]\,g'(u)\,\Phi(u)^{\,n}.
\]
Taking $g(u)=(u+1)^{a}$ and $\Phi(u)=(u+1)^2$ gives
\[
[z^{n}]\,\mathcal C(z)^{a}
=\frac{1}{n}\,[u^{\,n-1}]\,a(u+1)^{a-1}\,(u+1)^{2n}
=\frac{a}{n}\,[u^{\,n-1}]\,(u+1)^{2n+a-1}.
\]
Extracting the binomial coefficient yields
\[
[z^{n}]\,\mathcal C(z)^{a}
=\frac{a}{n}\binom{2n+a-1}{\,n-1\,} = \frac{a}{2n+a}\binom{2n+a}{n}.
\]

Let $o_r = 2n_r, r=1, 2, \ldots, a.$ Suppose each $\tau^{\mathbf o}_r$ is pure. It is known that $F(\tau^{\mathbf o}_r)=\gamma^{n_r}C_{n_r}.$ Therefore, 
\begin{equation*}
\begin{split}
\sum_{\substack{o_1+\cdots+o_a=p-a\\ o_r\in 2\mathbb N_0}} 
\ \prod_{r=1}^a F\big(\tau^{\mathbf o}_r\big) & = \sum_{\substack{n_1+\cdots+n_a=(p-a)/2\\ n_r\in \mathbb N_0}}
\ \prod_{r=1}^a F\big(\tau^{\mathbf o}_r\big)\\
& = \gamma^{(p-a)/2} \cdot \left[z^{(p-a)/2}\right]\mathcal C(z)^a.
\end{split}
\end{equation*}
On the other hand, suppose each $\tau^{\mathbf o}_r$ is alternating. It follows that $F(\tau^{\mathbf o}_r)=C_n.$ Hence,
\begin{equation*}
\begin{split}
\sum_{\substack{o_1+\cdots+o_a=p-a\\ o_r\in 2\mathbb N_0}} 
\ \prod_{r=1}^a F\big(\tau^{\mathbf o}_r\big) = \left[z^{(p-a)/2}\right]\mathcal C(z)^a.
\end{split}
\end{equation*}
In summary, \eqref{eq:outer-compress-F} follows by \eqref{eq:Lagrange-inversion}. 
\end{proof}

\begin{proof}[Proof of Theorem~\ref{thm:spoke-arc-closed}]
Let $X$ be the normalized complex Gaussian elliptic matrix given in \eqref{eq:normalized-elliptic}, and let $W_1, W_2$ be given in \eqref{eq:W}. By Corollary \ref{cor:spoke-arc}, we have
\begin{equation}\label{eq:proof-thm-spoke-arc-closed-1}
\lim_{N\to\infty}\kappa_2(W_1,W_2)
=\sum_{a=1}^{\min\{p,q\}}
 \sum_{\mathbf o,\boldsymbol{\iota}}
 \sum_{[(U,V)]\in \mathcal O_{a,\mathbf o,\boldsymbol\iota}}
\gamma^{\,s_{\mathrm{sp}}(a;\tau;U,V)}
\prod_{r=1}^{a} F\big(\tau^{\mathbf o}_r\big)\ \prod_{r=1}^{a} F\big(\tau^{\boldsymbol\iota}_r\big).
\end{equation}

Hence, the limiting covariance is written as a finite sum over the number of
spokes $a$, arc-lengths $\mathbf o, \boldsymbol{\iota}$ on the two circles and $\mathbb Z_a$-orbits of spoke endpoint lists $(U,V)$.
For each fixed $a$, the product of single-arc weights on the outer (resp.\ inner) circle is
given by Lemma~\ref{lem:arc-compress-unified}. So it remains to deal with the spoke factor
$\sum_{[(U,V)]}\gamma^{s_{\mathrm{sp}}}$ for each family.

\begin{itemize}
\item \textit{Pure powers vs.\ pure powers.}
In this case, $\tau(t) =1$ for any $t \in [p+q],$ i.e., all circle type words are all-$1$. 
Hence, for given $a$, $\tau^{\mathbf o}_r$ and $\tau^{\boldsymbol\iota}_r$ are pure. Moreover, every spoke is same-type, so we have $s_{\mathrm{sp}}=a$ for
any $\mathbb Z_a$-orbits $(U,V)$. Therefore, by Lemma \ref{lem:arc-compress-unified}, we obtain
\begin{equation*}
\begin{split}
\text{The RHS of} \; \eqref{eq:proof-thm-spoke-arc-closed-1} & =\sum_{a=1}^{\min\{p,q\}}
 \sum_{\mathbf o,\boldsymbol{\iota}}
 \sum_{[(U,V)]\in \mathcal O_{a,\mathbf o,\boldsymbol\iota}}
\gamma^{a}
\prod_{r=1}^{a} F\big(\tau^{\mathbf o}_r\big)\ \prod_{r=1}^{a} F\big(\tau^{\boldsymbol\iota}_r\big)\\
& =\sum_{a=1}^{\min\{p,q\}} \frac{pq}{a}\gamma^{a}
 \sum_{\mathbf o,\boldsymbol{\iota}} \prod_{r=1}^{a} F\big(\tau^{\mathbf o}_r\big)\ \prod_{r=1}^{a} F\big(\tau^{\boldsymbol\iota}_r\big)\\
 & = \sum_{a=1}^{\min\{p,q\}} \frac{pq}{a}\,\gamma^{\frac{p+q}{2}}\, \mathrm{FC} \left(a,\tfrac{p-a}{2} \right) \cdot
  \mathrm{FC}\left(a,\tfrac{q-a}{2}\right)\\
 &=\sum_{a=1}^{\min\{p,q\}} a  \gamma^{\frac{p+q}{2}}\,\binom{p}{\tfrac{p-a}{2}} \binom{q}{\tfrac{q-a}{2}}.
\end{split}
\end{equation*}
where we use the fact $|\mathcal O_{a,\mathbf o,\boldsymbol\iota}| = \frac{pq}{a}.$

\item \textit{Pure power vs.\ adjoint power.}
In this case, $\tau(t) =1, t \in [p]$ and $\tau(t) =\ast, t \in p+[q],$ i.e., the outer circle is all-$1$ and the inner all-$\ast$. Note that $\tau^{\mathbf o}_r$ and $\tau^{\boldsymbol\iota}_r$ are also pure for given $a$. 
Since every spoke is mixed-type, so $s_{\mathrm{sp}}=0$ for any spoke labeling. Thus in this case, 
\begin{equation*}
\begin{split}
\text{The RHS of} \; \eqref{eq:proof-thm-spoke-arc-closed-1} & =\sum_{a=1}^{\min\{p,q\}}
 \sum_{\mathbf o,\boldsymbol{\iota}}
 \sum_{[(U,V)]\in \mathcal O_{a,\mathbf o,\boldsymbol\iota}} \prod_{r=1}^{a} F\big(\tau^{\mathbf o}_r\big)\ \prod_{r=1}^{a} F\big(\tau^{\boldsymbol\iota}_r\big)\\
& =\sum_{a=1}^{\min\{p,q\}} \frac{pq}{a}
 \sum_{\mathbf o,\boldsymbol{\iota}} \prod_{r=1}^{a} F\big(\tau^{\mathbf o}_r\big)\ \prod_{r=1}^{a} F\big(\tau^{\boldsymbol\iota}_r\big)\\
 & = \sum_{a=1}^{\min\{p,q\}} \frac{pq}{a}\,\gamma^{\frac{p+q}{2}-a}\, \mathrm{FC} \left(a,\tfrac{p-a}{2} \right) \cdot
  \mathrm{FC}\left(a,\tfrac{q-a}{2}\right)\\
 &=\sum_{a=1}^{\min\{p,q\}} a \gamma^{\frac{p+q}{2}-a}\,
  \binom{p}{\tfrac{p-a}{2}}
  \binom{q}{\tfrac{q-a}{2}}.
\end{split}
\end{equation*}
which gives \eqref{eq:FCpa}.

\item \textit{Alternating vs.\ alternating.}
In this case, the outer and inner circle words are $(XX^\ast)^p$ and $(XX^\ast)^q$, respectively. So types on every single arc on the inner circle (resp. outer circle) are alternating, and contribute only the Catalan factors $\mathrm{FC}(a,\,p-\tfrac{a}{2})$ (resp. $\mathrm{FC}(a,\,q-\tfrac{a}{2})$) in the sum. Note that two circles have $2p$ and $2q$ points, respectively. 

For the spokes: because each arc length on both circles is even, the spoke endpoints alternate in type ($1,\ast,1,\ast,\dots$) along each circle. Define the \emph{unit rotation} $T$ on the inner endpoint list $U=(u_1,\dots,u_a)$ by
\(
T(U)\coloneqq(u_2,\dots,u_a,u_1).
\)
Consider the map
\[
\Phi:\ \mathcal O_{a,\mathbf o,\boldsymbol\iota}\to \mathcal O_{a,\mathbf o,\boldsymbol\iota},\qquad
[(U,V)]\longmapsto[(T(U),V)].
\]
\emph{Well-definedness and bijectivity.}
Let $R$ be the generator of the simultaneous cyclic relabeling on spoke indices:
$(U,V)\sim (R\!\cdot\!U,R\!\cdot\!V)$. Since $T$ and $R$ are both cyclic shifts on the inner list,
they commute: $T(R^k\!\cdot\!U)=R^k\!\cdot\!T(U)$ for all $k$. Hence for representatives in the same orbit,
\[
[(U,V)]=[(R^k\!\cdot\!U,R^k\!\cdot\!V)]
\ \Longrightarrow\
\Phi([(U,V)])=\Phi([(R^k\!\cdot\!U,R^k\!\cdot\!V)]),
\]
so $\Phi$ is well-defined on $\mathcal O_{a,\mathbf o,\boldsymbol\iota}$. The inverse map is $[(U,V)]\mapsto[(T^{-1}(U),V)]$; therefore $\Phi$ is a bijection.

Because types alternate on each circle, we have
$s_{\mathrm{sp}}(a;\tau;U,V)\in\{0,a\}$. Moreover, $\tau(u_{r+1})\neq\tau(u_r)$ for all $r$, hence
\[
s_{\mathrm{sp}}(a;\tau;T(U),V)=
\begin{cases}
0, & s_{\mathrm{sp}}(a;\tau;U,V)=a,\\[2pt]
a, & s_{\mathrm{sp}}(a;\tau;U,V)=0.
\end{cases}
\]
Thus $\Phi$ flips the spoke class: it maps every same-type orbit to a mixed-type orbit and vice versa. Since $\Phi$ is a bijection on these two classes, exactly half of the orbits contribute weight $\gamma^a$ and half contribute weight $1$.  In addition, we require $a$ to be even (since the numbers of arc points $2p-a$ and $2q-a$ are even), which we encode by the factor
$\frac{1+(-1)^a}{2}$. Therefore,
\[
\sum_{[(U,V)]}\gamma^{\,s_{\mathrm{sp}}}
\;=\;
\frac{2p\cdot 2q}{a}\cdot\frac{\gamma^a+1}{2}\cdot\frac{1+(-1)^a}{2}
\;=\;
\frac{2pq}{a}\,(\gamma^a+1)\,\frac{1+(-1)^a}{2}.
\]

Putting these pieces together,
\begin{equation*}
\begin{split}
&\text{The RHS of} \; \eqref{eq:proof-thm-spoke-arc-closed-1} \\
& \quad =\sum_{a=1}^{2\min\{p,q\}}
 \frac{2pq}{a}\cdot \frac{1+(-1)^a}{2}\,(\gamma^a+1) \sum_{\mathbf o,\boldsymbol{\iota}}
\prod_{r=1}^{a} F\big(\tau^{\mathbf o}_r\big)\ \prod_{r=1}^{a} F\big(\tau^{\boldsymbol\iota}_r\big)\\
&\quad =\sum_{a=1}^{2\min\{p,q\}} \frac{pq}{a}
 [1+(-1)^a] \,(\gamma^{a}+1) \ \mathrm{FC}\left(a,\,p-\tfrac{a}{2}\right) \cdot
   \mathrm{FC}\left(a,\,q-\tfrac{a}{2}\right)\\
 &\quad =\sum_{a=1}^{\min\{p,q\}} a(\gamma^{2a}+1)
\binom{2p}{p-a}\binom{2q}{q-a}.
\end{split}
\end{equation*}
which proves \eqref{eq:FCall}.
\end{itemize}
\end{proof}

\begin{proof}[Proof of Theorem~\ref{thm:spoke-arc-closed-real}]
The proof is similar to the complex case. Let $X$ be the normalized real Gaussian elliptic matrix given in \eqref{eq:normalized-elliptic}, and let $W_1, W_2$ be given in \eqref{eq:W}. By Corollary \ref{cor:real-spoke-arc}, we have
\begin{equation}
\begin{aligned}
		&\lim_{N\to\infty}\kappa_2(W_1,W_2)\\
&=\left[(\text{complex})+\sum_{a=1}^{\min\{p,q\}}
 \sum_{\mathbf o,\boldsymbol{\iota}}
 \sum_{[(U,V)]\in \mathcal O_{a,\mathbf o,\boldsymbol\iota}}
\gamma^{a-s_{sp}(a;\tau;U,V)}
\prod_{r=1}^{a}F\left( \tau^{\mathbf o}_r \right)\prod_{r=1}^{a}F\left(\tau^{\boldsymbol\iota}_r\right)\right].
\end{aligned}
\end{equation}
Here ``(complex)'' denotes the right-hand side of~\eqref{eq:spoke-arc-complex}. For fixed $a$, the products of single-arc weights are given by Lemma~\ref{lem:arc-compress-unified}. So we only need to consider the contribution of spokes.

For every configuration $(U, V),$ the spoke weight is the sum of the two leading contributions:
\[
\text{all pairs cross}:\ \gamma^{\,s_{\mathrm{sp}}}
\quad\text{and}\quad
\text{spokes straight, arcs cross}:\ \gamma^{\,a-s_{\mathrm{sp}}}.
\]
Therefore, the spoke factor over orbits is
\[
\sum_{[(U,V)]}\big(\gamma^{\,s_{\mathrm{sp}}}+\gamma^{\,a-s_{\mathrm{sp}}}\big).
\]
We now evaluate this in the three typical families; all other steps are identical to the complex case, and we omit the details. 

\begin{itemize}
\item \textit{Pure powers vs.\ pure powers.} The type words on both circles are all-$1$. Then $s_{\mathrm{sp}}=a$ for any $(U,V)$, so
\[
\sum_{[(U,V)]}\big(\gamma^{\,s_{\mathrm{sp}}}+\gamma^{\,a-s_{\mathrm{sp}}}\big)
=\frac{pq}{a}\,(\gamma^{a}+1).
\]

\item \textit{Pure power vs.\ adjoint power.}
The type word is all-$1$ on outer circle, and all-$\ast$ on inner circle. Then $s_{\mathrm{sp}}=0$ for any $(U,V)$, so the spoke
sum is again
\(
\frac{pq}{a}\,(\gamma^{a}+1).
\)

\item 
\textit{Alternating vs.\ alternating.}
Only even $a$ are admissible, and along the spoke endpoints the types are also alternate. Recall that in this case we have $s_{\mathrm{sp}}(a;\tau;U,V)\in\{0,a\}$. Therefore, 
\(
\gamma^{\,s_{\mathrm{sp}}}+\gamma^{\,a-s_{\mathrm{sp}}}=\gamma^{a}+1.
\)
Hence the spoke factor over $\Z_a$-orbits of spoke labelings is
\[
\sum_{[(U,V)]}\big(\gamma^{s_{\mathrm{sp}}}+\gamma^{a-s_{\mathrm{sp}}}\big)
=\frac{4pq}{a}\cdot \frac{1+(-1)^a}{2}(\gamma^{a}+1).
\]
\end{itemize}
\end{proof}

\subsection{Proof of Theorem~\ref{thm:elliptic-r2}}\label{subsec:second-free}
Here we present a self-contained, diagrammatic proof based on the spoke-arc decomposition, which also makes the transpose (outer-reversal) contribution in the real setting transparent. We only give the proof of the real case, and the same strategy holds for the complex case.

To isolate the combinatorial effect from centered clusters and cyclically alternating sequences, we first rewrite each trace as a product of centered single-color clusters and assume alternation along each trace. This effect removes cluster means and, in the large-$N$ limit, rules out contributions from pair partitions that keep an entire cluster internally paired. 

Write the two trace words as trace products of centered single-color clusters
\[
W_1=\Tr(A_1\cdots A_p),\qquad
W_2=\Tr(B_1\cdots B_q),
\]
with
\[
A_i\coloneqq Y_i-\alpha_i I,\quad
B_j\coloneqq Y_{p+j}-\alpha_{p+j} I,\quad 
\alpha_i\coloneqq\E\tr(Y_i)\ \ (1\le i\le p+q),
\]
where each $Y_i\in\mathrm{alg}\langle X^{(c_i)},(X^{(c_i)})^{*}\rangle$ is a single-color monomial of real elliptic matrices.
Assume both words are cyclically alternating in color, i.e.\
$c_i\neq c_{i+1}$ for $i\in[p]$ (indices mod $p$) along $W_1$ and 
$c_{p+j}\neq c_{p+j+1}$ for $j\in[q]$ (indices mod $q$) along $W_2.$

Define the letter intervals
\[
I_1:=[\,|Y_1|\,],\quad
I_2:=|Y_1|+[\,|Y_2|\,],\ \dots,\ 
I_i:=\big(|Y_1|+\cdots+|Y_{i-1}|\big)+[\,|Y_i|\,],
\]
so that $I_i$ records the positions occupied by $Y_i$ in the concatenation $Y_1Y_2\cdots Y_{p+q}$. By a slight abuse of terminology, we will not distinguish between a cluster and its letter interval: we may refer to $I_k$ (on the circle) and also to the associated monomial $A_k$ simply as the cluster.
\begin{lemma}\label{lem:center-alter-effect}
With the above notations, one has 
\[
\lim_{N\to\infty}\kappa_2\big(\Tr(A_1\cdots A_p),\ \Tr(B_1\cdots B_q)\big)
=\sum_{\pi\in S^\star}\gamma^{\,s(\pi;\tau)}+\gamma^{\,s(\pi;\tau)+a(\pi)-2\,s_{\mathrm{sp}}(\pi;\tau)},
\]
where $S^\star \subset NC_2^{c}\big(\cup_{i\in[p]} I_i,\ \cup_{j\in p+[q]} I_j\,\big)$ denotes the set of color-respecting non-crossing annular pair partitions in which no cluster is paired entirely internally; equivalently, for every cluster $I_k$ at least one of its letters is paired with a letter in $I_\ell$ for some $\ell\neq k$.
\end{lemma}

\begin{proof}
Expanding $A_i=Y_i-\alpha_i I$ and $B_j=Y_{p+j}-\alpha_{p+j}I$ and using the bilinearity of $\kappa_2$
gives the identity
\begin{equation}\label{eq:cluster-PIE-M}
\kappa_2\big(W_1,W_2\big)
=\sum_{M\subseteq[p+q]}
(-1)^{|M|}\Big(\prod_{t\in M}\alpha_t\Big)\,
\kappa_2\Big(\Tr Y_{[p]\setminus M}, \Tr Y_{p+[q]\setminus M}\Big).
\end{equation}
Here the cyclic products follow the trace orders of $W_1$ and $W_2$, respectively:
\[
Y_{[p]\setminus M}\coloneqq\prod_{i\in [p]\setminus M}^{\circlearrowleft}\! Y_i,\quad
Y_{p+[q]\setminus M}\coloneqq\prod_{i\in (p+[q])\setminus M}^{\circlearrowleft}\! Y_i,
\]
with the conventions $Y_{\varnothing}=I$ and $(\cdot)^{\circlearrowleft}$ denoting cyclic products. Note that $[p]\setminus M$ is a shorthand for $[p]\setminus(M\cap[p])$, and similarly $(p+[q])\setminus M=(p+[q])\setminus\big(M\cap(p+[q])\big)$.

For $M\subseteq[p+q]$, define the inner and outer letter sets
\[
P_M\coloneqq\bigcup_{i\in [p]\setminus M}\! I_i,\quad
Q_M\coloneqq\bigcup_{j\in (p+[q])\setminus M}\! I_j.
\]
Let $\tau_M$ be the global type word on $P_M\cup Q_M$ induced by
$Y_{[p]\setminus M}$ and $Y_{p+[q]\setminus M}$.
Then, by Proposition~\ref{prop:semi-closed-cov-real} we have
\begin{align}\notag
&\kappa_2\Big(\Tr Y_{[p]\setminus M}, \Tr Y_{p+[q]\setminus M}\Big)
\\ \notag
&\qquad\quad =\sum_{\pi\in NC_2^c(P_M,Q_M)}
\left(
\gamma^{\,s(\pi;\tau_M)}
+\gamma^{\,s(\pi;\tau_M)+a(\pi)-2s_{\mathrm{sp}}(\pi;\tau_M)}
\right)+O(N^{-1}).
\end{align}
Here $NC_2^c(P_M,Q_M)$ denotes the color-respecting non-crossing annular pair partitions on the two circles with sets $P_M$ and $Q_M$.
For $\pi\in NC_2^c(P_M,Q_M)$, $a(\pi)$ is the number of spokes of $\pi$,
$s(\pi;\tau_M)$ counts all same-type pairs in $\pi$ and
$s_{\mathrm{sp}}(\pi;\tau_M)$ is the number of same-type spokes in $\pi$.

Let $S=NC_2^c(\cup_{i\in[p]}I_i,\cup_{j\in p+[q]}I_j)$ and $E_k=\{\pi\in S:I_k\ \text{is paired internally}\}$.
For brevity, we use the shorthand $\Omega(\pi;\tau)$ for the real elliptic weight:
\[
\Omega(\pi;\tau):=\gamma^{\,s(\pi;\tau)}+\gamma^{\,s(\pi;\tau)+a(\pi)-2s_{\mathrm{sp}}(\pi;\tau)}.
\]
Then for any $M\subseteq[p+q]$,
\[
\sum_{\pi\in \bigcap_{k\in M} E_k} \Omega(\pi;\tau)
=\Big(\prod_{k\in M}\alpha_k\Big)
\sum_{\pi'\in NC_2^{c}(P_M,Q_M)} \Omega(\pi';\tau_M),
\]
and, by inclusion-exclusion principle (see, e.g. \cite[Section III.7]{FlajoletSedgewick2009}),
\[
\sum_{M\subseteq[p+q]}(-1)^{|M|}\!\sum_{\pi\in \bigcap_{k\in M}E_k}\Omega(\pi;\tau)
=\sum_{\pi\in S\setminus\bigcup_{k=1}^{p+q}E_k}\Omega(\pi;\tau).
\]
Let
\[
S^\star\coloneqq S\setminus\bigcup_{k=1}^{p+q}E_k.
\]
$S^\star$ is precisely the subset of $S$ in which no cluster is paired entirely internally. Hence, only those color-respecting non-crossing annular pair partitions in which \emph{at least one letter in every cluster is paired with a letter from another cluster} contribute to the limiting covariance.
\end{proof}

\begin{proof}[Proof of Theorem~\ref{thm:elliptic-r2}]
Recall that the inner cycle is on $[p]$ and the outer cycle is on $p+[q]$. $\rho=(1, 2, \ldots, p)\,(p+1, p+2, \ldots,p+q)\in S_{p+q}$ is the trace permutation. Let $\tau: [p+q] \to \{1,\ast\}$ be a global type word. 

\smallskip
\noindent\underline{\emph{First-order freeness.}}
By Adhikari--Bose \cite[Theorem~2]{adhikari2019brown}, a family of independent real elliptic matrices are asymptotically free. Thus the first-order part is settled.

\smallskip
\noindent\underline{\emph{Second-order limit distribution.}}
As an immediate consequence of the spoke-arc formulas, elliptic families admit a second-order limit distribution. Indeed, the first-order functional $\varphi_1$ is given by the single-arc limit $F_c(\cdot)$ (see Remark~\ref{rmk:multi-color}), while the second-order functional $\varphi_2$ is given explicitly by Corollary \ref{cor:real-spoke-arc}. 
So it remains to prove the vanishing of higher-order cumulants.

Fix $m\ge3$ trace words $W_1,\dots,W_m$ (finite length, arbitrary colors). Let
$\ell_r:=|W_r|$ and $E:=\sum_{r=1}^m\ell_r$ be the total number of matrix factors.
By Wick's formula and independence across colors, only \emph{color-respecting} pair partitions $\pi\in\mathcal P_2^{\mathrm{c}}(E)$ survive (each pair joins entries of the same color). For each pair $\{x,y\}\in\pi$ (of color $c$), we use \eqref{eq:elliptic-entry-cov-real}. Expanding the product over pairs yields a sum over \emph{assignments} $\varepsilon\in\{c,s\}^{E/2}$ choosing cross/straight on each pair:
\begin{equation}\notag
\mathbb{E}\Big[\prod_{r=1}^m W_r\Big]
=\sum_{\pi\in\mathcal P_2^{\mathrm{c}}(E)}\sum_{ \varepsilon}
\Big(\tfrac{1}{N}\Big)^{E/2}
\sum_{\mathbf i} \prod_{\{x,y\}\in\pi}
\delta^{\varepsilon}_{xy} \prod_{\{x,y\}\in\pi} w_{\varepsilon(\{x,y\})}\big(\tau(x),\tau(y)\big),
\end{equation}
where $\delta^{c}_{xy}:=\delta_{i_x i_{\rho(y)}}\delta_{i_{\rho(x)} i_y}$ and
$\delta^{s}_{xy}:=\delta_{i_x i_{y}}\delta_{i_{\rho(x)} i_{\rho(y)}}$.

For fixed $(\pi,\varepsilon)$ the product
$\prod_{\{x,y\}\in\pi}\delta^{\varepsilon}_{xy}$ identifies indices along the cycles of an \emph{index permutation} $\theta=\theta(\pi; \varepsilon)$ (see \eqref{def:index-permut}). Thus
\[
\sum_{\mathbf i} \prod_{\{x,y\}\in\pi}\delta^{\varepsilon}_{xy} = N^{\#(\theta(\pi; \varepsilon))}.
\]
Hence the contribution of $(\pi,\varepsilon)$ equals
\begin{equation}\label{eq:N-power-index}
\mathrm{Contr}(\pi,\varepsilon)=N^{\#(\theta(\pi; \varepsilon))-E/2} \prod_{\{x,y\}\in\pi} w_{\varepsilon(\{x,y\})}\big(\tau(x),\tau(y)\big).
\end{equation}

Passing from moments to cumulants uses the classical \emph{moment-cumulant formula} on the lattice of set partitions $\mathcal P(m)$ (see, e.g., \cite[Section~1]{mingo2017free}):
\begin{equation}\label{eq:MS-mc}
\kappa_m(W_1,\dots,W_m)
=\sum_{\alpha\in\mathcal P(m)} \mu_{\mathcal P}(\alpha,\mathbf 1_m)
\prod_{V\in\alpha} \mathbb E\Big[\prod_{r\in V} W_r\Big],
\end{equation}
where $\mathbf 1_m$ is the one-block partition and $\mu_{\mathcal P}$ is the Möbius function of the poset $\mathcal P(m)$. Let $J_r\subset [E]$ be the set of positions belonging to $W_r$, i.e., with $s_r:=1+\sum_{t<r}\ell_t$ and $e_r:=\sum_{t\le r}\ell_t$, we have $J_r=\{s_r,s_r+1,\dots,e_r\}\subset [E]$. If $\alpha$ is a partition of $[m]$, define $\tilde\alpha\in\mathcal P(E)$ to be the partition of $[E]$ such that whenever $r$ and $s$ lie in the same block of $\alpha$, then $J_r$ and $J_s$ lie in the same block of $\tilde\alpha$. Given $\pi\in \mathcal P(E)$, we define $\hat{\pi}$ be the partition on $[m]$ such that $r$ and $s$ are in the same block of $\hat{\pi}$ if there is a block of $\pi$ contains $J_r$ and $J_s$.

Next we regroup the cluster moment associated to a specific partition $\alpha\in\mathcal P(m)$. Write
$J_V:=\bigsqcup_{r\in V}J_r\subset[E]$, $E_V:=\sum_{r\in V}\ell_r$ for $V\in\alpha$.
If $\pi\le \tilde\alpha$ (i.e., no pair of $\pi$ links different
$J_V$'s), then $\pi$ decomposes uniquely as a disjoint union
\[
\pi=\bigsqcup_{V\in\alpha}\pi_V\quad\text{with}\ \ \pi_V\in\mathcal P_2(J_V),
\]
and the assignment $\varepsilon$ decomposes accordingly as $\varepsilon=(\varepsilon_V)_{V\in\alpha}$.
Since the delta functions do not couple indices across different $J_V$'s, we have
\[
\#(\theta(\pi;\varepsilon))=\sum_{V\in\alpha}\#(\theta(\pi_V;\varepsilon_V)),
\qquad
E=\sum_{V\in\alpha}E_V,
\]
and the type-weights factor clusterwise. Hence
\[
\mathrm{Contr}(\pi,\varepsilon)
=\prod_{V\in\alpha}
\Bigg(
N^{\#(\theta(\pi_V;\varepsilon_V))-E_V/2}
\!\!\!\prod_{\{x,y\}\in\pi_V} w_{\varepsilon_V(\{x,y\})}\big(\tau(x),\tau(y)\big)
\Bigg).
\]
Summing over all $(\pi,\varepsilon)$ with $\pi\le\tilde\alpha$ thus factorizes:
\[
\sum_{\substack{(\pi,\varepsilon)\\ \pi\le \tilde\alpha}}\mathrm{Contr}(\pi,\varepsilon)
=\prod_{V\in\alpha}
 \sum_{\pi_V\in\mathcal P_2(J_V)}\ \sum_{\varepsilon_V\in\{c,s\}^{E_V/2}}
\cdots.
\]
The inner double sum is precisely the Wick expansion of the cluster moment on $V$, i.e.
\[
\sum_{\pi_V,\varepsilon_V}\cdots=\ \mathbb E\Big[\prod_{r\in V} W_r\Big].
\]
Therefore
\begin{equation}\label{eq:cluster-moment}
\sum_{\substack{(\pi,\varepsilon)\\ \pi\le \tilde\alpha}}\mathrm{Contr}(\pi,\varepsilon)
=\prod_{V\in\alpha} \mathbb E\Big[\prod_{r\in V} W_r\Big].
\end{equation}

Insert the Wick expansion of each cluster moment in \eqref{eq:MS-mc}. Using \eqref{eq:cluster-moment} we obtain
\begin{align}\notag
\kappa_m(W_1,\dots,W_m)
&=\sum_{\alpha\in\mathcal P(m)}\mu_{\mathcal P}(\alpha,\mathbf 1_m)
   \sum_{\substack{(\pi,\varepsilon)\\ \pi\le \tilde \alpha}}
   \mathrm{Contr}(\pi,\varepsilon)
   \\ \notag
&=\sum_{(\pi,\varepsilon)} \mathrm{Contr}(\pi,\varepsilon)
   \sum_{\substack{\alpha\in\mathcal P(m)\\ \alpha\ge \hat\pi}}
   \mu_{\mathcal P}(\alpha,\mathbf 1_m).
   \notag
\end{align}
Note that only the terms with $\hat\pi=\mathbf 1_m$ (i.e., connected) contribute, hence
\[
\kappa_m(W_1,\dots,W_m)
=\sum_{(\pi,\varepsilon):\hat\pi=\mathbf 1_m}
\mathrm{Contr}(\pi,\varepsilon).
\]
   
For such $(\pi,\varepsilon)$ we use the following Euler characteristic formula.
\[
\#(\theta(\pi;\varepsilon))-\frac{E}{2}\ \le\ 2-2g-m,\quad g\in\mathbb N_0,
\]
with equality iff all arcs use the cross constraint, and all spokes are of the same type
(either all cross or all straight). This follows from the index identifications proved in Section~\ref{Semi-closed-lim-cov} (arc-straight merges two index cycles; mixing cross/straight on spokes also merges cycles).

Therefore each connected term in \eqref{eq:N-power-index} satisfies
\[
N^{\,\#(\theta(\pi;\varepsilon))-E/2}\ \le\ N^{\,2-2g-m}\ \le\ N^{\,2-m},
\]
which implies 
\[
\lim_{N\to\infty}\kappa_m(W_1,\dots,W_m)=0
\]
for $m\ge3.$

\medskip 

In the rest of the proof, we will verify Equation \eqref{eq:second-order-free-defin}. 

\medskip
\noindent\textit{Admissible patterns for $\pi \in S^\ast$.}
Recall that 
$
S^\star=S\setminus\bigcup_{k=1}^{p+q}E_k
$
is the subset of $S$ such that at least one letter in every cluster is paired with a letter from another cluster. Lemma~\ref{lem:center-alter-effect} shows that, in the large-$N$ limit, only those $\pi\in S^\star$ contribute to the covariance. In what follows we carry out a case-by-case analysis of $S^\star$ subject to the additional constraints, isolating the configurations that survive in the large-$N$ limit.

First, all arc pairs in $\pi\in S^\ast$ lie within a single cluster. Indeed, assume there is an arc pair joining letters from $I_i$ and $I_j$ with $i<j$, and choose such an pair with minimal number of internal clusters $j-i$. If $j=i+1$, alternation forces $c(I_i)\neq c(I_{i+1})$, contradicting color-respecting. If $j-i\ge2$, consider the inner region strictly between $I_i$ and $I_j$. By non-crossing, letters in this region cannot connect outside across the chosen pair. Since we work in $S^\star$, no cluster in the region is entirely self-paired, hence within the region there exists a cross-cluster arc pair that connects adjacent clusters. however, this is also impossible by alternation. The outer circle can be treated similarly. Therefore, any $\pi \in S^\star$ contains no  cross-cluster arc pairs.

Since there is at least one cross-cluster pair in any $\pi\in S^\ast$ and all arc pairs in any $\pi\in S^\ast$ lies within a single cluster, cross-cluster pairs must be spokes, and thus each cluster carries at least one spoke. In words, $\pi\in S^\star$ contributes if and only if \emph{every cluster carries at least one spoke and all remaining pairs are within-cluster arc pairs.}

\medskip
\noindent\textit{Coherence of spokes (from the letter-level to the cluster-level).}
Fix an inner circle cluster $A_k$ with interval $I_k$ and hence $I_k$ carries at least one spoke endpoint. Let $J_k$ denote the set of spoke endpoints on the outer circle that are connected to letters of $I_k$ by spokes.

\emph{Claim 1: all spoke endpoints in $J_k$ lie inside a single outer cluster; i.e., there exists a unique map $\sigma:[p]\to p+[q]$ such that for every $k\in[p]$, $J_k\subseteq I_{\sigma(k)}$.} Recall that $I_k$ carries at least one spoke endpoint. Suppose $J_k$ meets two distinct outer clusters. Since $J_k\neq\varnothing$, pick $u,v\in J_k$ such that $J_k\subseteq [u,v]_\rho \coloneqq \{u,\rho(u),\ldots,v\}$ on the outer circle. Let $I_\alpha$ and $I_\beta$ denote the outer clusters containing $u$ and $v$, respectively; that is, $u\in I_\alpha$ and $v\in I_\beta$. We consider two cases:

\begin{itemize}
\item \textit{Case 1: $I_\alpha$ and $I_\beta$ are adjacent.}
By alternation, $c(I_\alpha)\neq c(I_\beta)$. Since a spoke must be color-respecting, every spoke issued from $I_k$ can land only in an outer cluster of the same color $c(I_k)$. Hence at most one of $I_\alpha$ and $I_\beta$ can receive a spoke endpoint from $I_k$, which contradicts our choice of $u,v\in J_k$ with $u\in I_\alpha$ and $v\in I_\beta$.

\item \textit{Case 2: $I_\alpha$ and $I_\beta$ are not adjacent and $I_\alpha\neq I_\beta$.} Since $u,v\in J_k$, color-respecting spokes force $c(I_\alpha)=c(I_\beta)=c(I_k)$. For any outer cluster $I_m$ between $I_\alpha$ and $I_\beta$ along the outer circle, $S^\star$ ensures that $I_m$ contains at least a spoke endpoint. By non-crossing, any spoke issued from $I_m$ cannot cross the spokes from $I_k$ into $I_\alpha$ and $I_\beta$, hence its inner endpoint must lie on the subinterval of $I_k$ between those two inner endpoints. Therefore that spoke also connects to $I_k$, so $c(I_m)=c(I_k)$ for any $I_m$ between $I_\alpha$ and $I_\beta$. However, by color alternation on the outer circle, adjacent clusters have different colors; in particular $c(I_{\alpha+1})\neq c(I_\alpha)=c(I_k)$, a contradiction.
\end{itemize}

Both cases are impossible, so $J_k$ cannot meet two distinct outer clusters (i.e., $I_\alpha=I_\beta$), and there exists a unique map $\sigma:[p]\to p+[q]$ such that for every $k\in[p]$, $J_k\subseteq I_{\sigma(k)}$.

\emph{Claim 2. The map $\sigma: [p]\rightarrow p+[q]$ is a bijection and thus $p=q$. Moreover, $\sigma$ reverses cyclic order, i.e., $\sigma(k+1)=\sigma(k)-1\pmod p.$}
The bijection of $\sigma$ follows from the outer-inner symmetric version of Claim~1: defining $\zeta:p+[q]\to[p]$ analogously yields $\zeta\circ\sigma=\mathrm{id}_{[p]}$ and $\sigma\circ\zeta=\mathrm{id}_{p+[q]}$, hence $\sigma$ is a bijection and $p=q$. Under the non-crossing property, the spokes must appear in the same geometric order on the annulus.
Because we read the inner circle clockwise while the outer circle is read counterclockwise, ``same geometric order'' translates into ``reversed cyclic order'' in the indexings, hence $\sigma(k+1)=\sigma(k)-1\ (\mathrm{mod}\ p)$.

Combining Claim~1 and 2 with the preceding discussion, we obtain the following characterization of the patterns that contribute in the limit: \emph{Contributing patterns $\pi\in S^\star$ are exactly those for which there exists a unique bijection $\sigma:[p]\to p+[p]$ that reverses the cyclic order and such that, for every $k$, all spokes from $I_k$ terminate in the unique outer cluster $I_{\sigma(k)}$. Within each cluster, after cutting at the spoke endpoints, each resulting subinterval carries a (possibly empty) non-crossing pair partition.}

Figure~\ref{fig:contributing-pattern} depicts a contributing pattern $\pi\in S^\ast$ in two views: a global schematic and a local diagram for a matched pair of clusters. In the global schematic we show only the spokes connecting matched clusters, suppressing the non-crossing pair partitions on the arcs; the local diagram reinstates these pair partitions, with the remainder of the annulus omitted.

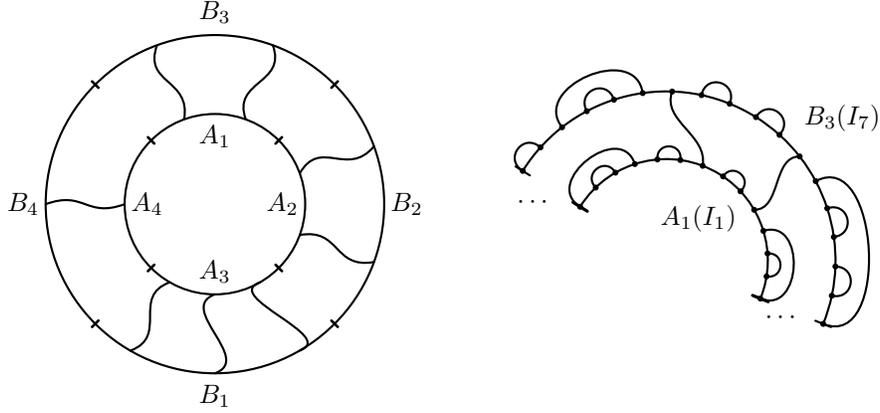
\begin{figure}[H]
\centering
\begin{tikzpicture}[line cap=round,line join=round,thick,scale=.75]

\begin{scope}
  \def\R{3}   
  \def\r{1.6}   

  \draw (0,0) circle (\R);
  \draw (0,0) circle (\r);

\def\tick{0.08} 
\foreach \m in {45,135,225,-45}{
  \draw[line width=1pt] (\m:{\r-\tick}) -- (\m:{\r+\tick});
  \draw[line width=1pt] (\m:{\R-\tick}) -- (\m:{\R+\tick});
}

    \node at (180:\r-0.40) {$A_4$};
    \node at ( 90:\r-0.40) {$A_1$};
    \node at ( 0:\r-0.40) {$A_2$};
    \node at (-90:\r-0.40) {$A_3$};

    \node at (180:\R+0.40) {$B_4$};
    \node at ( 90:\R+0.40) {$B_3$};
    \node at ( 0:\R+0.40) {$B_2$};
    \node at (-90:\R+0.40) {$B_1$};

  \draw (70:\r) to[out=110, in=-70, looseness=1.35]   (70:\R);  
  \draw (110:\r) to[out=80, in=-110, looseness=1.35]   (110:\R);    
  \draw ( -20:\r) to[out= 10, in=190, looseness=1.35]  (-20:\R);   
  \draw ( 20:\r) to[out=60, in=-130, looseness=1.35] (20:\R);     
  \draw (-60:\r) to[out=200, in= 20, looseness=1.35]   (-60:\R);  
  \draw (-90:\r) to[out=200, in= 20, looseness=1.35]   (-90:\R);
  \draw (-120:\r) to[out=200, in= 20, looseness=1.35]   (-120:\R);
  \draw (-180:\r) to[out=200, in= 20, looseness=1.35]   (-180:\R);
\end{scope}

\begin{scope}[xshift=8cm, yshift=-1cm]
  \def\R{3.0}        
  \def\r{1.8}        
  \def\astart{148}    
  \def\aend{-22}     
  \def\bend{-8}
  \draw (0,0) ++(\astart:\R) arc (\astart:\aend:\R);
  \draw (0,0) ++(\astart:\r)  arc (\astart:\aend:\r);

  \foreach \ang in {148,138,128,118,108,98,88,78,68,58,48,38,28,18,8,-2,-12,-22}
    {\fill (\ang:\R) circle (1.5pt);}
  \foreach \ang in {148,134,120,108,95,82,69,56,43,30,17,4,-9,-22}
    {\fill (\ang:\r) circle (1.5pt);}
  
\def\tick{0.15} 
\foreach \m in {148,-22}{
  \draw[line width=1.2pt] (\m:{\r-\tick}) -- (\m:{\r+\tick});
  \draw[line width=0.8pt] (\m:{\R-\tick}) -- (\m:{\R+\tick});
}  

  \draw (30:\r) to[out=0,in=180,looseness=1] (38:\R);    
  \draw (69:\r) to[out=80,in=-110,looseness=1] (88:\R);
  \draw (148:\r) to[out=148,in=108,looseness=1.7] (108:\r); 
  \draw (134:\r) to[out=134,in=120,looseness=2] (120:\r);
  \draw (95:\r) to[out=95,in=82,looseness=2] (82:\r);
  \draw (56:\r) to[out=56,in=43,looseness=2] (43:\r);
  \draw (17:\r) to[out=17,in=-22,looseness=1.7] (-22:\r);
  \draw (4:\r) to[out=4,in=-9,looseness=2] (-9:\r);
  \draw (-22:\R) to[out=-22,in=28,looseness=1.3] (28:\R); 
  \draw (-2:\R) to[out=-2,in=-12,looseness=2] (-12:\R);
  \draw (18:\R) to[out=18,in=8,looseness=2] (8:\R);
  \draw (48:\R) to[out=48,in=58,looseness=2] (58:\R);
  \draw (68:\R) to[out=68,in=78,looseness=2] (78:\R);
  \draw (98:\R) to[out=98,in=128,looseness=1.5] (128:\R);
  \draw (108:\R) to[out=108,in=118,looseness=2] (118:\R);
  \draw (138:\R) to[out=138,in=148,looseness=2] (148:\R);
  
  \node[anchor=west] at (48:\R+0.40) {$B_3 (I_7)$};
  \node[anchor=west] at (110:\r-1.) {$A_1 (I_1)$};

  \node at (156:2.55) {$\cdots$};
  \node at (-26:2.3) {$\cdots$};
\end{scope}
\end{tikzpicture}
\caption{Left: a global schematic with inner clusters $A_1,\dots,A_4$ and outer clusters $B_1,\dots,B_4$; spokes realize the order-reversing matching $\sigma$ (here $\sigma(1)=3$, $\sigma(2)=2$, $\sigma(3)=1$, $\sigma(4)=4$). Right: a local diagram of the matched clusters $A_1$ and $B_3$ connected by two spokes; after cutting at the spoke endpoints, each subinterval on the arcs carries a non-crossing pair partition. Here $|I_1|=14$ and $|I_{7}|=18$.}\label{fig:contributing-pattern}
\end{figure}

\noindent\textit{Passage to cluster-level weights.}
Fix $k\in [p], \sigma$ and a contributing pattern \(\pi\in S^\ast\). Let us consider the local subdiagram of $\pi$ that involves only the connected intervals $I_k$ and $I_{\sigma(k)}$, which correspond to the clusters $A_k$ and $B_{\sigma(k)}$ in our setup. We then glue the endpoints of $I_k$ and $I_{\sigma(k)}$ pairwise to obtain a single circle (see Figure~\ref{fig:glue} for this local gluing). 

Reading the outer interval in the reversed cyclic order yields a single circle with  $|I_k|+|I_{\sigma(k)}|$ marked points, on which the restriction of $\pi$ induces a non-crossing pair partition \(\pi^\sigma_k\in NC_2(I_k\cup I_{\sigma(k)})\). The induced type word is $\tau^\sigma_k\coloneqq\tau\big|_{I_k}\times \tau\big|_{I_{\sigma(k)}}$. Since $S^\star$ forces both $I_k$ and $I_{\sigma(k)}$ to carry at least one spoke, the resulting non-crossing partition $\pi_k^\sigma$ has at least one pair connecting $I_k$ and $I_{\sigma(k)}$.

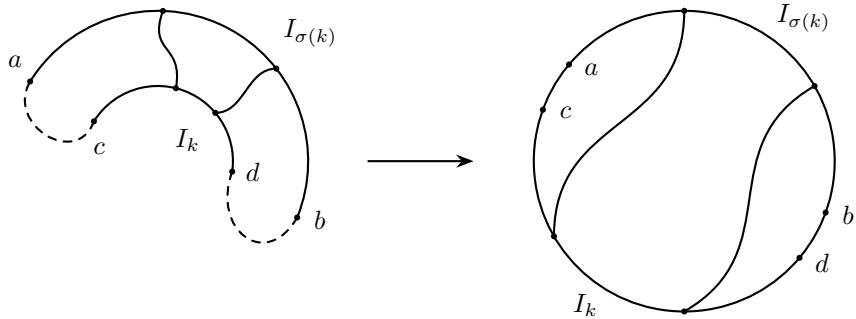
\begin{figure}[H]
\centering
\begin{tikzpicture}[line cap=round,line join=round,thick,scale=1]
\usetikzlibrary{arrows.meta}
  \def\R{2.0}        
  \def\r{1}        
  \def\astart{148}    
  \def\aend{-22}     
  \def\bend{-8}
  \draw (0,0) ++(\astart:\R) arc (\astart:\aend:\R);
  \draw (0,0) ++(\astart:\r)  arc (\astart:\bend:\r);

  \foreach \ang in {148,88,38,-22}
    {\fill (\ang:\R) circle (1.15pt);}
  \foreach \ang in {148,76,40,-8}
    {\fill (\ang:\r) circle (1.15pt);}

  \draw (40:\r) to[out=0,in=180,looseness=1] (38:\R);    
  \draw (76:\r) to[out=80,in=-110,looseness=1.5] (88:\R);
  \draw[dashed] (-8:\r) to[out=-110,in=-120,looseness=2] (-22:\R);
  \draw[dashed] (148:\r) to[out=-110,in=-120,looseness=1.7] (148:\R);

  \node[anchor=west] at (48:\R+0.3) {$I_{\sigma(k)}$};
  \node[anchor=west] at (70:\r-0.7) {$I_k$};
  \node[anchor=west] at (148:\R+0.50) {$a$};
  \node[anchor=west] at (170:\r-0) {$c$};
  \node[anchor=west] at (-22:\R+0.10) {$b$};
  \node[anchor=west] at (-8:\R-0.95) {$d$};
  
  \coordinate (Ledge) at (\R+0.8,0);

  \begin{scope}[xshift=7cm] 
    \def\RR{2} 
    \coordinate (Redge) at (-\RR-0.8,0);

    \draw (0,0) circle (\RR);
    \foreach \ang in {210,270,30,90,140,160,-20,-40}
      {\fill (\ang:\RR) circle (1.15pt);}
    \draw (210:\RR) to[out=90,in=-90,looseness=1.2] (90:\RR);
    \draw (270:\RR) to[out=30,in=210,looseness=1.2] (30:\RR);
    \node[anchor=west] at (140:\RR-0.1) {$a$};
    \node[anchor=west] at (160:\RR-0.1) {$c$};
    \node[anchor=west] at (-20:\RR+0.10) {$b$};
    \node[anchor=west] at (-40:\RR+0.10) {$d$};
    \node[anchor=west] at (60:\RR+0.2) {$I_{\sigma(k)}$};
    \node[anchor=west] at (230:\RR+0.5) {$I_k$};
    
  \end{scope}

  \draw[->,>=Stealth,thick] (Ledge) -- (Redge);
\end{tikzpicture}
\caption{Local gluing of connected clusters $I_k$ and $I_{\sigma(k)}$ (two spokes shown, non-crossing pair partitions on arcs suppressed). The endpoints $a,b$ of $I_{\sigma(k)}$ and $c,d$ of $I_k$ are glued pairwise ($a\leftrightarrow c$, $b\leftrightarrow d$) to form a single circle.}
\label{fig:glue}
\end{figure}

Globally, any fixed \(\pi\in S^\ast\) decomposes into \(p\) local pair partitions on \(I_k\cup I_{\sigma(k)}\): 
\[
\pi = \pi_1^{\sigma} \times \cdots \times \pi_p^{\sigma},\qquad 
\pi_k^{\sigma}\in NC_2\big(I_k\cup I_{\sigma(k)}\big),\ \ k=1,\dots,p,
\]
where \(\sigma:[p]\to p+[q]\) is the unique bijection given in the above {\it Claim 2}. It reverses the cyclic order determined by \(\pi\)
and each \(\pi_k^\sigma\) contains at least one pair connecting \(I_k\) and \(I_{\sigma(k)}\).
In particular, we have
\[
s(\pi;\tau)=\sum_{k=1}^{p}s\big(\pi_k^{\sigma};\tau_k^{\sigma}\big),\qquad
a(\pi)=\sum_{k=1}^{p}a\big(\pi_k^{\sigma}\big),\qquad
s_{\mathrm{sp}}(\pi;\tau)=\sum_{k=1}^{p}s_{\mathrm{sp}}\big(\pi_k^{\sigma};\tau_k^{\sigma}\big).
\]

Now we are ready to prove Equation \eqref{eq:second-order-free-defin}. Recalling Lemma~\ref{lem:center-alter-effect} we have,
\begin{equation}\label{eq:zero-term}
\lim_{N\to\infty}\kappa_2\big(\Tr(A_1\cdots A_p),\ \Tr(B_1\cdots B_q)\big)
=\sum_{\pi\in S^\star}\gamma^{\,s(\pi;\tau)}+\gamma^{\,s(\pi;\tau)+a(\pi)-2\,s_{\mathrm{sp}}(\pi;\tau)}.
\end{equation}
For the first term in the summation, we invoke the first-moment identity \eqref{id-arcweight-moment} and use the same argument in Lemma~\ref{lem:center-alter-effect}. For any centered clusters $A_k,B_{\sigma(k)}$, we obtain
\[
\lim_{N\to\infty}\E\tr\big(A_k B_{\sigma(k)}\big)=\sum_{\substack{\pi\in NC_2(I_k\cup I_{\sigma(k)})\\\text{at least one cross-cluster pair}}}\gamma^{\,s(\pi;\tau_k^\sigma)}=\sum_{\pi_k^\sigma}\gamma^{\,s(\pi_k^\sigma;\tau_k^\sigma)}.
\]
Then, summing over all possible $\sigma$ yields
\begin{equation}\label{eq:frist-term}
\sum_{\pi\in S^\star}\gamma^{\,s(\pi;\tau)}=\sum_\sigma\prod_{k=1}^p\sum_{\pi_k^\sigma}\gamma^{\,s(\pi_k^{\sigma};\tau_k^{\sigma})}=\sum_{i=1}^p\prod_{k=1}^p\lim_{N\to\infty}\E\tr\big(A_k B_{k-i}\big).
\end{equation}

For the second term (mixed assignment), we impose the same ``at least one cross-cluster pair'' restriction. By Remark~\ref{spoke-straight=outer-reversed-transposed}, the mixed assignment (for the real ensemble) is indexwise identical to the complex all-cross channel after replacing the outer product by its transpose read in reversed order; accordingly we arrange the outer factors as $B_q^T,\dots,B_1^T$ clockwise. This arrangement induces the order-preserving bijection $\eta$ between inner and outer clusters, i.e., $\eta(k+1)=\eta(k)+1 (\mathrm{mod}\ p)$. The same decomposition gives
\[
\sum_{\pi\in S^\star}\gamma^{\,s(\pi;\tau)+a(\pi)-2\,s_{\mathrm{sp}}(\pi;\tau)}=\sum_\eta\prod_{k=1}^p\sum_{\pi_k^\eta}\gamma^{\,s(\pi_k^{\eta};\tau_k^{\eta})+a(\pi_k^{\eta})-2s_{\mathrm{sp}}(\pi_k^{\eta};\tau_k^{\eta})},
\]
Applying the first-moment identity to the transpose term yields, for centered clusters,
\[
\sum_{\pi_k^\eta}\gamma^{\,s(\pi_k^{\eta};\tau_k^{\eta})+a(\pi_k^{\eta})-2s_{\mathrm{sp}}(\pi_k^{\eta};\tau_k^{\eta})}=\lim_{N\rightarrow\infty}\E\tr \big(A_k B^T_{\eta(k)}\big).
\]
Here we note that the first-moment limit is the same for both real and complex elliptic ensembles.

Hence
\begin{equation}\label{eq:second-term}
\sum_{\pi\in S^\star}\gamma^{\,s(\pi;\tau)+a(\pi)-2\,s_{\mathrm{sp}}(\pi;\tau)}=\sum_\eta\prod_{k=1}^p\sum_{\pi_k^\eta}\gamma^{\,s(\pi_k^{\eta};\tau_k^{\eta})}=\sum_{i=1}^p\prod_{k=1}^p\lim_{N\rightarrow\infty}\E\tr \big(A_k B^T_{k+i}\big).
\end{equation}

Combining \eqref{eq:zero-term}, \eqref{eq:frist-term}, and \eqref{eq:second-term}, we deduce \eqref{eq:second-order-free-defin}, which completes our proof.
\end{proof}


\bigskip

{\bf Acknowledgment.} We thank James A. Mingo, Dylan Gawlak, and Yong Jiao for valuable discussions and comments. We are partially supported by NSFC 12031004, 12571147 and Provincial Natural Science Foundation of Hunan
(grant number 2024JJ1010). L. Wu further acknowledge the support of the Institut Henri Poincaré (UAR 839 CNRS-Sorbonne Université) and of LabEx CARMIN (ANR-10-LABX-59-01).

\bibliographystyle{acm}
\bibliography{matrix}

\begin{thebibliography}{10}

\bibitem{adhikari2019brown}
{\sc Adhikari, K., and Bose, A.}
\newblock Brown measure and asymptotic freeness of elliptic and related matrices.
\newblock {\em Random Matrices Theory Appl. 8}, 02 (2019), 1950007.

\bibitem{Bordenave-Chafai-circular}
{\sc Bordenave, C., and Chafa\"{\i}, D.}
\newblock Around the circular law.
\newblock {\em Probab. Surv. 9\/} (2012), 1--89.

\bibitem{Collins2007}
{\sc Collins, B., Mingo, J.~A., \'Sniady, P., and Speicher, R.}
\newblock Second order freeness and fluctuations of random matrices. {III}. {H}igher order freeness and free cumulants.
\newblock {\em Doc. Math. 12\/} (2007), 1--70.

\bibitem{Crisanti1987}
{\sc Crisanti, A., and Sompolinsky, H.}
\newblock Dynamics of spin systems with randomly asymmetric bonds: Langevin dynamics and a spherical model.
\newblock {\em Phys. Rev. A 36\/} (1987), 4922--4939.

\bibitem{FlajoletSedgewick2009}
{\sc Flajolet, P., and Sedgewick, R.}
\newblock {\em Analytic Combinatorics}.
\newblock Cambridge University Press, 2009.

\bibitem{Fyodorov2016}
{\sc Fyodorov, Y.~V., and Khoruzhenko, B.~A.}
\newblock Nonlinear analogue of the {M}ay-{W}igner instability transition.
\newblock {\em Proc. Natl. Acad. Sci. USA 113}, 25 (2016), 6827--6832.

\bibitem{Fyodorov1997}
{\sc Fyodorov, Y.~V., Khoruzhenko, B.~A., and Sommers, H.-J.}
\newblock Almost-{H}ermitian random matrices: eigenvalue density in the complex plane.
\newblock {\em Phys. Lett. A 226}, 1 (1997), 46--52.

\bibitem{Girko1986}
{\sc Girko, V.~L.}
\newblock Elliptic law.
\newblock {\em Theory Probab. Appl. 30}, 4 (1986), 677--690.

\bibitem{Girko1995}
{\sc Girko, V.~L.}
\newblock The elliptic law: ten years later \uppercase {I}.
\newblock {\em Random Oper. Stoch. Equ. 3}, 3 (1995), 257--302.

\bibitem{Girko19952}
{\sc Girko, V.~L.}
\newblock The elliptic law: ten years later \uppercase {II}.
\newblock {\em Random Oper. Stoch. Equ. 3}, 4 (1995), 377--398.

\bibitem{Girko1997}
{\sc Girko, V.~L.}
\newblock Strong elliptic law.
\newblock {\em Random Oper. Stoch. Equ. 5}, 3 (1997), 269--306.

\bibitem{Girko2006}
{\sc Girko, V.~L.}
\newblock The strong elliptic law, twenty years later, \uppercase {I}.
\newblock {\em Random Oper. Stoch. Equ. 14}, 1 (2006).

\bibitem{Girko20062}
{\sc Girko, V.~L.}
\newblock The strong elliptical galactic law. sand clock density, twenty years later, \uppercase {II}.
\newblock {\em Random Oper. Stoch. Equ. 14}, 2 (2006), 157--208.

\bibitem{Lando2004-dc}
{\sc Lando, S.~K., and Zvonkin, A.~K.}
\newblock {\em Graphs on Surfaces and Their Applications}, vol.~141 of {\em Encyclopaedia of mathematical sciences}.
\newblock Springer-Verlag, Berlin, 2004.

\bibitem{Marti2018}
{\sc Mart\'{\i}, D., Brunel, N., and Ostojic, S.}
\newblock Correlations between synapses in pairs of neurons slow down dynamics in randomly connected neural networks.
\newblock {\em Phys. Rev. E 97\/} (2018), 062314.

\bibitem{MingoNica2004}
{\sc Mingo, J.~A., and Nica, A.}
\newblock Annular noncrossing permutations and partitions, and second-order asymptotics for random matrices.
\newblock {\em Int. Math. Res. Not. IMRN}, 28 (2004), 1413--1460.

\bibitem{MingoPopa2013}
{\sc Mingo, J.~A., and Popa, M.}
\newblock Real second order freeness and {H}aar orthogonal matrices.
\newblock {\em J. Math. Phys. 54}, 5 (2013), 051701, 35.

\bibitem{Mingo2007}
{\sc Mingo, J.~A., \'Sniady, P., and Speicher, R.}
\newblock Second order freeness and fluctuations of random matrices. {II}. {U}nitary random matrices.
\newblock {\em Adv. Math. 209}, 1 (2007), 212--240.

\bibitem{MingoSpeicher2006}
{\sc Mingo, J.~A., and Speicher, R.}
\newblock Second order freeness and fluctuations of random matrices. {I}. {G}aussian and {W}ishart matrices and cyclic {F}ock spaces.
\newblock {\em J. Funct. Anal. 235}, 1 (2006), 226--270.

\bibitem{mingo2017free}
{\sc Mingo, J.~A., and Speicher, R.}
\newblock {\em Free probability and random matrices}, vol.~35.
\newblock Springer, 2017.

\bibitem{Naumov12}
{\sc Naumov, A.}
\newblock Elliptic law for real random matrices.
\newblock {\em arXiv:1201.1639\/} (2012).

\bibitem{NguyenORourke2015}
{\sc Nguyen, H.~H., and O'Rourke, S.}
\newblock The elliptic law.
\newblock {\em Int. Math. Res. Not. IMRN}, 17 (2015), 7620--7689.

\bibitem{nica2006lectures}
{\sc Nica, A., and Speicher, R.}
\newblock {\em Lectures on the combinatorics of free probability}, vol.~13.
\newblock Cambridge University Press, 2006.

\bibitem{Poley2024}
{\sc Poley, L., Galla, T., and Baron, J.~W.}
\newblock Eigenvalue spectra of finely structured random matrices.
\newblock {\em Phys. Rev. E 109\/} (Jun 2024), 064301.

\bibitem{Raney1960}
{\sc Raney, G.~N.}
\newblock Functional composition patterns and power series reversion.
\newblock {\em Trans. Amer. Math. Soc. 94}, 3 (1960), 441--451.

\bibitem{Redelmeier2014}
{\sc Redelmeier, C. E.~I.}
\newblock Real second-order freeness and the asymptotic real second-order freeness of several real matrix models.
\newblock {\em Int. Math. Res. Not. IMRN}, 12 (2014), 3353--3395.

\bibitem{Schomerus2017}
{\sc Schomerus, H.}
\newblock Random matrix approaches to open quantum systems.
\newblock {\em Stochastic Processes and Random Matrices: Lecture Notes of the Les Houches Summer School 2015\/} (2017), 409--473.

\bibitem{Sommers1988}
{\sc Sommers, H.-J., Crisanti, A., Sompolinsky, H., and Stein, Y.~J.}
\newblock Spectrum of large random asymmetric matrices.
\newblock {\em Phys. Rev. Lett. 60\/} (1988), 1895--1898.

\bibitem{SuryaWarnke2023}
{\sc Surya, E., and Warnke, L.}
\newblock Lagrange inversion formula by induction.
\newblock {\em Am. Math. Mon. 130}, 10 (2023), 944--948.

\bibitem{TaoVu2010}
{\sc Tao, T., and Vu, V.}
\newblock Random matrices: universality of {ESD}s and the circular law.
\newblock {\em Ann. Probab. 38}, 5 (2010), 2023--2065.
\newblock With an appendix by Manjunath Krishnapur.

\bibitem{tutte1962census}
{\sc Tutte, W.~T.}
\newblock A census of slicings.
\newblock {\em Canad. J. Math. 14\/} (1962), 708--722.

\bibitem{Voiculescu1991}
{\sc Voiculescu, D.}
\newblock Limit laws for random matrices and free products.
\newblock {\em Invent. Math. 104}, 1 (1991), 201--220.

\bibitem{Voiculescu1998}
{\sc Voiculescu, D.}
\newblock A strengthened asymptotic freeness result for random matrices with applications to free entropy.
\newblock {\em Int. Math. Res. Not. IMRN}, 1 (1998), 41--63.

\end{thebibliography}

\end{document}